\UseRawInputEncoding
\documentclass[11pt, leqno]{amsart}
\usepackage{amsmath, amssymb, latexsym}
\usepackage{mathrsfs}
\usepackage{amsmath}
\usepackage{enumerate}
\usepackage{color}
\usepackage[colorlinks=true, pdfstartview=FitV, linkcolor=blue,%
citecolor=blue, urlcolor=blue]{hyperref}
\usepackage[all, knot]{xy}
\usepackage{chngcntr}
\counterwithout{figure}{section}
\xyoption{arc}
\usepackage{tikz}
\usepackage{mathdots}
\usepackage{amsmath}
\allowdisplaybreaks

\numberwithin{equation}{section}

\newtheorem{theorem}{Theorem}[section]
\newtheorem{corollary}[theorem]{Corollary}
\newtheorem{lemma}[theorem]{Lemma}
\newtheorem{proposition}[theorem]{Proposition}

\theoremstyle{definition}
\newtheorem{definition}[theorem]{Definition}

\newcommand{\mbf}{\mathbf}
\newcommand{\mbb}{\mathbb}
\newcommand{\mcal}{\mathcal}

\title[The Extension groups of one $\rm K$-type representations of graded Hecke algebra]
{ The Extension groups of one $\rm K$-type representations of graded Hecke algebra}

\author[Quanyong Chen]{Quanyong Chen}
\address{The University of Hong Kong,
	Hong Kong, China}
\email{chenquanyong@hrbeu.edu.cn}

\thanks{ }

\begin{document}
\begin{abstract}
We study the extension groups between one $\rm K$-type representations of the graded Hecke algebra $\mathbb H$.
A one $\rm K$-type representation is a Hermitian representation of $\mbb H$ whose restriction to the Weyl group is irreducible.
Using a Koszul-type resolution of $\mathbb H$-modules, we reduce the computation of $\operatorname{Ext}_{\mathbb H}^j(X,Y)$ to a combinatorial count, yielding explicit formulas for the classical types and a computer computation for the exceptional types.
\end{abstract}
	
	\maketitle

	\section{Introduction}
	\subsection{Background}
	
	The graded Hecke algebra $\mbb H=\mbb H(\Pi,k)$, introduced by Lusztig \cite{Lu88, Lu89, Lu95} and fundamental results can be found in \cite{Op95, HO97, OS10}. It serves as a bridge between the serves as a bridge between the representation theory of Weyl groups and that of $p$-adic groups \cite{Bo76,KL87}.
	The algebra carries a natural $\ast$-structure.
	Following Barbasch and Moy \cite{BM99}, a Hermitian representation $\pi$ of $\mbb H$ is called a one $\rm K$-type representation if its restriction to $W$ is irreducible.
    One $\rm K$-type representations play a central role in the Barbasch-Moy approach to the unitary dual of
    $p$-adic groups \cite{BM89,BM99}.

    The study of Ext-groups of one $\rm K$-type representations arises from the continuous cohomology of unitary representations \cite{BW00}.
	In the real case, the basic problem is to determine when the continuous cohomology $H_{\rm ct}^\ast(G,H_\pi\otimes F)$ of a unitary representation $H_\pi$ of a real Lie group $G$ with coefficients in a finite-dimensional representation $F$ does not vanish.
	Vogan and Zuckerman \cite{VZ84} reduced this analytic problem to an algebraic one by identifying the continuous cohomology with the $(\mathfrak{g},K)$-cohomology of the underlying Harish-Chandra module,
	\begin{equation*}
		H^\ast_{\rm ct}(G,H_\pi)\cong H^\ast_{(\mathfrak{g},K)}(H_\pi^{K\text{-fin}}),
	\end{equation*}
	and the latter can in turn be computed as $\operatorname{Ext}^i_{(\mathfrak{g},K)}(\mbb C_{\rm triv},H_\pi^{K\text{-fin}})$.
	In the $p$-adic case, the Koszul-type resolution \cite{Ch16} takes the place of the $(\mathfrak{g},K)$ Cartan--Eilenberg resolution, and the corresponding invariants are the extension groups of one $\rm K$-type representations of the graded Hecke algebra, whose computation is the subject of this paper.
	The corresponding computation for $p$-adic group was carried out by Casselman \cite{Ca74}.
	 For a semisimple $p$-adic group $G$,
	\begin{align*}
	\operatorname{Ext}^i_G(\mbb C_{\rm triv},\operatorname{St})\cong
	\begin{cases}
		\mbb C & i=\operatorname{rank}G,\\
		0 & i\neq \operatorname{rank}G,
	\end{cases}
	\end{align*}
	where $\operatorname{St}$ is the irreducible subquotient of $C^\infty(B\backslash G)$ attached to the trivial character of a Borel subgroup $B$.
	The graded Hecke algebra version of this result is obtained from our formulas by taking $X$ and $Y$ to be the one $\rm K$-type
	representations attached to the trivial and the sign characters of $W$, which recovers the same answer in degree $j={\rm rank}G$. In this sense the results of this paper generalize the theorem of Casselman.
	
	Orlik \cite{Or05} computed the extension groups of a closely related class of representations, the generalized Steinberg representations, that is, the irreducible composition factors of the principal series at the character $\chi(v)=\langle v,\rho^\vee\rangle$ of $S(V)$, where $\rho^\vee$ denotes the half-sum of the positive coroots. Since a one $\rm K$-type module is annihilated by $V$ while a generalized Steinberg representation has nonzero central character, the two classes of modules are disjoint, and Orlik's results complement rather than overlap with ours.
	
\subsection{Main results}
The main results are collected in the following theorem, which provides an efficient combinatorial algorithm for computing the extension groups between arbitrary one $\rm K$-type representations, expressed in terms of Young diagrams and the coefficients of a polynomial.

In type $A$, the generating function $\prod_{i=1}^{s}(1+q^{2i-1})$
enumerates partitions of $m$ into distinct odd parts chosen from
$\{1,3,\dots,2s-1\}$, which are in bijection with self-conjugate establishes a direct connThese partitions are in bijection with the self-conjugate partitions contained in the $s\times s$ square, so the Ext groups in type $A$ are directly tied to the combinatorics of self-conjugate partitions.

\begin{theorem}\label{thm:typeA}
	{\rm (Proposition \ref{theorem 1.1.1}, Theorem \ref{theorem 1.1.2} and Proposition \ref{theorem 1.1.3}.)} Let $X,Y$ be the one $\rm K$-type representations of graded Hecke algebra $\mbb H_{A_{n-1}}$, corresponding to the rectangular partitions $\lambda=(k^d)$ and $\mu=(k'^{d'})$. Then $\operatorname{Ext}_{\mbb H_{A_{n-1}}}^j(X,Y)=0$ for all $j$ unless $\{k,d\}=\{k',d'\}$. In the nontrivial case $\mu=(d^k)$, with $s=\min(d,k)$,
	$$\dim\operatorname{Ext}_{\mbb H_{A_{n-1}}}^j(X,Y)=[q^{n-j}]\prod_{i=1}^{s}(1+q^{2i-1}),$$
	while for self-extensions $\dim\operatorname{Ext}_{\mbb H_{A_{n-1}}}^j(X,X)=[q^{j}]\prod_{i=1}^{s}(1+q^{2i-1})$.
\end{theorem}

In type $B$, the main combinatorial tool is the Mackey decomposition for the double coset space of parabolic subgroups of $B_n$.
It expresses the dimension of ${\rm Ext}_{\mbb H_{B_n}}^j(X,Y)$ as the cardinality of the set  $\Pi_{\alpha,\beta}^{j}$, which is in bijection with the partitions contained in a certain rectangle.

\begin{theorem}\label{thm:typeBC}
{\rm (Proposition \ref{both one partition}, Proposition \ref{one-two bipartition}, Theorem \ref{bipartition multiplication of B}.)}
Let $X,Y$ be the one $\rm K$-type representations of the graded Hecke algebra $\mbb H_{B_n}$,  corresponding to the rectangular bipartitions $\alpha=(a^p)(b^q)$ and $\beta=(c^k)(d^l)$, respectively.  Then
$\operatorname{Ext}_{\mbb H_{B_n}}^j(X,Y)=0$ for all $j$ unless $\alpha, \beta$ satisfy one of the cases in \eqref{four conditions}. Then we have
$$ {\rm dim \ Ext}_{\mbb H_{B_n}}^j(X,Y)=
[q^j] q^{j_0} [W+H, H]_{q^2},$$
If one of the two bipartitions is empty, the Ext groups vanish except in a single degree, where the dimension is $1$.
\end{theorem}

In type $D$, the new phenomenon is that the restriction of a self-associate representation of $B_n$ to $D_n$ splits into two inequivalent irreducible components.
 The Ext groups in type $D_n$ are obtained by combining the type $B_n$ formulas with this restriction, as stated in the following theorem.

\begin{theorem}\label{thm:typeD}
	{\rm (Proposition \ref{theorem 1,3.1}, Proposition \ref{theorem 1.3.2}.)} Let $X, Y$ be the one $\rm K$-type representations of graded Hecke algebra.
	The formulas are obtained from those for type $B_n$ by restriction. For non-split representations, $\dim\operatorname{Ext}_{\mbb H}^j(X,Y)=\sharp\Pi^j_{\alpha,\beta}+\sharp\Pi^j_{\alpha,\overline{\beta}}$, where $\overline{\beta}=(d^l)(c^k)$.
	For split representations, only even degrees $j=2r$ occur, and the dimensions are expressed in terms of the difference characters $\Delta_\lambda=\chi_{\lambda^+}-\chi_{\lambda^-}$ and the parity sign $\theta(-1)^r$ with $\theta\in\{+,-\}$.
\end{theorem}

The exceptional types are treated by computer.
For types $E_6, E_7$ and $E_8$, the restriction to $W$ of a one $\rm K$-type representation is self-twisted exactly when its character vanishes on the simple reflection $s_1$, and otherwise twisting by the sign character changes the sign of this value.
For type $F_4$, the pair of values on the simple reflections $s_1$ and $s_3$ takes the place of this single value.
For type $G_2$, the four one  $\rm K$-type representations of dimension one occur for all parameters, and a fifth one of dimension two appears when $h=1$ or $h=3$.
For each exceptional type we list the irreducible representations  appearing in $\chi\otimes \wedge^iV$ that are the restrictions of one $\rm K$-type representations.
Since the differentials vanish, the Ext groups are read off directly from these tables.
	
\subsection{organization}
The paper is organized as follows. Section $2$
 collects the graded Hecke algebra setup, the realizations of the classical Weyl group representations, and Koszul-type resolution together with the vanishing of the differentials.
In Section $3$, the Ext groups in type $A$ are computation for type $B$ goes through the Mackey decomposition of the double coset space, and the resulting formula counts pairs of subpartitions satisfying rectangle-complement conditions.
In Section $5$, the formulas for type $D$ are obtained from those of type $B$ by restriction, with difference characters handling the split case.
In Section $6$, the exceptional types are treated by computer, and the Ext groups are read off from the decompositions of $\chi\otimes \wedge^iV$ over $\mbb C[W]$.
	
\noindent {\bf Acknowledgement.}
I am grateful to Kei Yuen Chan for helpful discussions and valuable comments.
I would like to thank ChatGPT,
especially for its assistance in generating the code for exceptional type calculations used in this paper.
 The generated code was subsequently reviewed and validated by the authors, who remain responsible for its correctness and integration.
\section{Preliminaries}

	\subsection{Graded Hecke algebra}
	Let $R$ be a reduced root system with a fixed set $\Pi$ of simple roots, which determines the set $R^+$ of positive roots.
     Let $W$ be the finite reflection group associated with $R$, and let $V_0$ be a real vector space containing $R$ (not necessarily spanned by $R$).
     For each $\alpha \in \Pi$, let $s_\alpha \in W$ be the simple reflection corresponding to $\alpha$, for which $\alpha$ lies in the $-1$-eigenspace of $s_\alpha$.
     For $\alpha \in R$, let $\alpha^\vee \in \operatorname{Hom}_{\mathbb R}(V_0, \mathbb R)$ be defined by
\begin{equation*}
s_\alpha(v)=v-\langle v,\alpha^\vee \rangle \alpha,
\end{equation*}
where $\langle v,\alpha^\vee \rangle=\alpha^\vee(v)$.
Let $R^\vee \subset \operatorname{Hom}_{\mathbb R}(V_0, \mathbb R)$ consisting of all such $\alpha^\vee$, and set $V_0^\vee = \operatorname{Hom}_{\mathbb R}(V_0, \mathbb R)$.
Extending scalars to $\mathbb C$, we define $V = \mathbb C \otimes_{\mathbb R} V_0$ and $V^\vee = \mathbb C \otimes_{\mathbb R} V_0^\vee$.
The tuple $(R, V, R^\vee, V^\vee, \Pi)$ is called a root datum.

Let $k: \Pi \rightarrow \mathbb C$ be a constant parameter function on each $W$-orbit of $\Pi$.
We set  $k_\alpha$ for $k(\alpha)$.

	\begin{definition}
		(\cite[\S 2]{Ch16}.) The graded Hecke algebra $\mbb H = \mbb H(\Pi, k)$ associated with a root datum $(R, V, R^\vee, V^\vee, \Pi)$ and a parameter function $k$ is a unital associative
		algebra over $\mbb C$, generated by the symbols $\{t_w : w \in W\}$ and $\{ f_v : v\in V \}$, subject to the following relations:
		\begin{enumerate}
			\item The map $w \mapsto t_w$ from $\mbb C[W] \rightarrow \mbb H$ is an algebra injection.
			\item The map $v \mapsto f_v$ from $S(V ) \rightarrow H$ is an algebra injection.\\
			We shall slightly abuse the notation and write $v$ in place of $f_v$ from now on.			
			\item  For $\alpha \in \Pi$ and $v \in V \subset S(V)$, the generators satisfy the cross-relation:
			\begin{align*}
			t_{s_\alpha}v-s_{\alpha}(v)t_{s_\alpha}=k_\alpha\langle v, \alpha^{\vee} \rangle.
			\end{align*}
		\end{enumerate}
	\end{definition}
When the Weyl group is specified, we write $\mbb H_W$ for the corresponding graded Hecke algebra, for example $\mbb H_{A_{n-1}}$, $\mbb H_{B_n}$ and $\mbb H_{D_n}$.

	For $\mbb H$-modules $X$ and $Y$, let ${\rm Ext}_{\mbb H}^i(X, Y)$ denote the $i$-th extension group in the category of $\mbb H$-modules.
	We also denote by ${\rm Res}_W$ the restriction functor from $\mbb H$-modules to $\mbb C[W]$-modules.

    Recall that from \cite{BM99}, the following proposition shows that $\mbb H$ carries a natural $\ast$-involution.
	\begin{proposition}\label{star-involution}
 The algebra $\mbb H$ is a $\ast$ algebra with the $\ast$ involution given by
		\begin{align*}
			w^\ast&=w^{-1}, \ \ t \in W, \\
			v^\ast&=-v+\sum_{\beta \in R^+}t_\beta\langle v,\check{\beta} \rangle k_\beta,  \ \ v \in V.
		\end{align*}
	\end{proposition}
	
	\begin{definition}
		The irreducible representations of $\mbb C[W]$ that can be extended as a Hermitian representation. We call such a representation a one $\rm K$-type representation.
	\end{definition}
	
	We use the standard bases and positive systems for each of the classical types $A, B, C, D$
	Denote the system by $(R,\Pi, \check{R}, \check{\Pi})$ and the positive root system by $R^+$. \\
	$\mbf {Type}$ $A_{n-1}$.
	\begin{equation}
		\begin{split}
			R=\check{R}&=\{\epsilon_i-\epsilon_j\}_{i\neq j \leq n},\\
			\Pi=\check{\Pi}&=\{\epsilon_i-\epsilon_{i+1}\}_{1\leq i\leq n-1},\\
			R^+&=\{\epsilon_i-\epsilon_j\}_{i<j\leq n}.
		\end{split}
	\end{equation}
	In practice, it is often more convenient to work with the lattice corresponding to \(\mathfrak{gl}(n)\) rather than \(\mathfrak{sl}(n)\); we shall do so without further comment.\\
	$\mbf {Type}$ $B_{n}$.
	\begin{equation}
		\begin{split}
			R&=\{\pm\epsilon_i \pm \epsilon_j\}_{i<j\leq n}\cup \{\pm \epsilon_i\}_{i\leq n} , \ \ \check{R}=\{\pm \epsilon_i\pm \epsilon_j \}_{i<j\leq n} \cup \{\pm 2\epsilon_i\}_{i\leq n},\\
			\Pi&=\{\epsilon_{i}-\epsilon_{i+1}, \epsilon_n\}_{1\leq i\leq n-1}, \ \ \check{\Pi}=\{ \epsilon_{i}-\epsilon_{i+1}, 2\epsilon_n\}_{1\leq i\leq n-1},\\
			R^+&=\{\epsilon_i\pm \epsilon_j\}_{i<j\leq n} \cup \{\epsilon_i\}_{i\leq n}.
		\end{split}
	\end{equation}
	$\mbf {Type}$ $C_{n}$.
	\begin{equation}
		\begin{split}
			R&=\{\pm\epsilon_i \pm \epsilon_j\}_{i<j\leq n}\cup \{\pm 2\epsilon_i\}_{i\leq n} , \ \ \check{R}=\{\pm \epsilon_i\pm \epsilon_j \}_{i<j\leq n} \cup \{\pm \epsilon_i\}_{i\leq n},\\
			\Pi&=\{\epsilon_{i}-\epsilon_{i+1}, 2\epsilon_n\}_{1\leq i\leq n-1}, \ \ \check{\Pi}=\{ \epsilon_{i}-\epsilon_{i+1}, \epsilon_n\}_{1\leq i\leq n-1},\\
			R^+&=\{\epsilon_i\pm \epsilon_j\}_{i<j\leq n} \cup \{2\epsilon_i\}_{i\leq n}.
		\end{split}
	\end{equation}
	$\mbf {Type}$ $D_{n}$.
	\begin{equation}
		\begin{split}
			R=\check{R}&=\{\pm\epsilon_i \pm \epsilon_j\}_{i<j\leq n}, \\
			\Pi=\check{\Pi}&=\{\epsilon_i-\epsilon_{i+1}, \epsilon_{n-1}+\epsilon_n\}, \\
			R^+&=\{\epsilon_i\pm \epsilon_j\}_{i<j\leq n} .
		\end{split}
	\end{equation}

    We follow \cite{Lu79} in describing the realizations of irreducible representations.
In general, for any representation $\sigma$, we attach a harmonic polynomial $\pi_\sigma$ with the property that the canonical $W$-action on $V_\sigma := \operatorname{span}\{w\cdot \pi_\sigma\}_{w\in W}$
yields a realization of $\sigma$.
Using the coordinates fixed earlier, for each subset $I \subseteq \{1,\dots,n\}$, we set
	\begin{equation}
		\begin{split}
			\pi_{A,I}&=\prod_{\substack{i<j \\ i,j \in I}}(\epsilon_i-\epsilon_j),\\
			\pi_{B,I}=\pi_{C,I}&=\prod_{\substack{i<j \\ i,j \in I}}(\epsilon_i^2-\epsilon_j^2) \cdot\prod \epsilon_i,\\
			\pi_{D,I}&=\prod_{\substack{i<j \\ i,j \in I}} (\epsilon_i^2-\epsilon_j^2).
		\end{split}
	\end{equation}
	When $I=\{1,\cdots,n\}$, the polynomial $\pi_{\ast, I}$ a?ords the $sgn$ representation of the corresponding Weyl group. \\
	$\mbf {Type}$ $A_{n-1}$.
	The representations of the Weyl group are indexed by partitions.
	Let $(n_1,\cdots,n_k)$ with $n_i\leq n_{i+1}$ be a partition of $n$.
 Fix a partition of $\{1,\cdots.n\}$ into subsets $I_i (i= 1,\cdots,k)$
	of sizes $n_i$, and set
	\begin{equation}
		\pi_{I_1,\cdots,I_k}=\prod_i \pi_{A, I_i}.
	\end{equation}
The span of these harmonic polynomial ranging over all partitions into  subsets of sizes $(n_1,\cdots,n_k)$ generates an irreducible representation $V_\sigma$ attached to $\sigma$.
	The action of $W$ is $s_w\pi_{I_1 ,\cdots,I_k}=\pi_{wI_1 ,\cdots,wI_k}$.\\
	$\mbf {Types}$ $B,C,D$.
	In types $B$ and $C$, irreducible representations are parametrized by a pair of partitions $(\sigma_L)(\sigma_R)$ with total size $n$.
For type $D$, the pairs $[\sigma_L][\sigma_R]$ and $[\sigma_R][\sigma_L]$ give the same representation, namely the restriction of the corresponding type $B_n$ or $C_n$ representation.
When $\sigma_L = \sigma_R$, this restriction decomposes into two inequivalent irreducible components, denoted $[\sigma]^+$ and $[\sigma]^{-}$.
Following the conventions for type $A$, partition $\{1,\dots,n\}$ into subsets $I_s$ and $J_r$ according to the parts of $\sigma_L$ and $\sigma_R$, respectively, and set
\begin{align*}
\pi_{I,J} := \prod_s \pi_{D,I_s} \cdot \prod_r \pi_{B,J_r}.
\end{align*}
The span of the $W$-orbits of $\pi_{I,J}$ generates the irreducible representation $V_\sigma$.

   As proved in \cite{BM99}, the following proposition characterizes the representations of $W(A)$ that lift to $\mbb H$.
	\begin{proposition}
  A representation $\sigma$ of $W(A)$ extends to a representation
		of $\mbb H$ if and only if $n_1=\cdots=n_k=d$.
	\end{proposition}
Also from \cite{BM99},	we now consider the cases of types  $B, C$ and $D$.
\begin{proposition}
If either $\sigma_L$ or $\sigma_R$ is empty, then $\sigma$ extends to a Hermitian representation of $\mbb H$.
\end{proposition}
	
Partition $\{1,\dots,n\}$ into subsets $J_1,\dots,J_s$ and $K_1,\dots,K_r$, where the sizes of the $J_i$'s and $K_j$'s are the parts of $\sigma_L$ and $\sigma_R$, respectively.
Define
\begin{align}\label{cc}
		\delta =
		\begin{cases}
			2k_{\epsilon_1} & \text{in type } B, \\
			k_{2\epsilon_1} & \text{in type } C, \\
			0 & \text{in type } D.
		\end{cases}
	\end{align}
	
	\begin{proposition}\label{conditionBCD}
		Let $\delta$ be as in \eqref{cc}.
		Assume that $\sigma_L \neq  \emptyset$ corresponds to the
		partition $m_1\leq \cdots \leq m_k$, and $\sigma_R \neq  \emptyset$ corresponds to the partition $n_1\leq \cdots\leq n_l$.
		Then $\sigma$ extends to an hermitian representation of $\mbb H$ if and only if all $m_i = d$, all
		$n_j= f$ and $k-d= l-f +\delta$.
	\end{proposition}
	
	\subsection{Koszul-type resolution on $\mbb H$-modules}
	Let $X$ be an $\mbb H$-module and let $U$ be a finite-dimensional $\mbb C[W]$-module.
    The tensor space $\mbb H \otimes_{\mbb C[W]} U$ is an $\mbb H$-module by the left multiplication on the first factor,  and  $\mbb H$ acts on the space ${\rm Hom}_W (\mbb H, U)$ by the right translation, explicitly that is for $f \in  {\rm Hom}_W (\mbb H, U)$, the action of $h\in \mbb H$ is given by:
	$$(h'. f)(h)=f(hh'),  \ \ \text{ for \ all } h \in \mbb H.$$
	
	\begin{lemma}(Frobenius reciprocity.)
	There are natural isomorphisms
		$${\rm Hom}_{\mbb H}(X, {\rm Hom}_W (\mbb H, U)) = {\rm Hom}_W ({\rm Res}_W X, U),$$
		and
		$${\rm Hom}_{\mbb H}(\mbb H \otimes_{\mbb C[W]} U, X) = {\rm Hom}_W (U, {\rm Res}_W X).$$
	\end{lemma}
	We recall from \cite[Lemma 3.2]{Ch16} that the induced module is projective and ${\rm Hom}_W (\mbb H, U)$
is injective.
	\begin{lemma}
	Let $U$ be a $\mbb C[W]$-module. Then $\mbb H \otimes _{\mbb C[W]} U$ is projective and ${\rm Hom}_W (\mbb H, U)$ is injective.
	\end{lemma}
	Let $X$ be an $\mbb H$-module. Define a sequence of $\mbb H$-maps $d_i$ as follows:
	\begin{equation}\label{Koszul-type-complex}
		\begin{split}
			0 \rightarrow \mbb H \otimes_{\mathbb{C}[W]} \big({\rm Res}_W X \otimes \wedge^n V \big)
			\xrightarrow{d_{n-1}}
			\cdots &\xrightarrow{d_i}
			\mbb H \otimes_{\mathbb{C}[W]} \big({\rm Res}_W X \otimes \wedge^i V \big)\\
			&\xrightarrow{d_{i-1}}
			\cdots \xrightarrow{d_0}
			\mbb H \otimes_{\mathbb{C}[W]} {\rm Res}_W X
			\xrightarrow{\epsilon} X \to 0
		\end{split}
	\end{equation}
	such that $\epsilon: \mbb H \otimes_{\mbb C[W]} {\rm Res}_W X \rightarrow X$ is given by
	$$\epsilon (h\otimes x)=h.x$$
	and for $i\geq 1$, $d_i : \mbb H \otimes_{\mbb C[W]} ({\rm Res}_W X \otimes  \wedge^{i+1}V ) \rightarrow \mbb H \otimes_{\mbb C[W]} ({\rm Res}_W X \otimes \wedge^i V )$ given by
	\begin{align*}
		d_i&(h\otimes (x\otimes v_1 \wedge \cdots \wedge v_{i+1}))\\
		&=\sum_{j=1}^{i+1}(-1)^{j+1}\times(hv_j\otimes x\otimes v_1\wedge \cdots\wedge \widehat{v_j}\wedge \cdots v_{i+1}-h\otimes v_j.x \otimes v_1\wedge \cdots \widehat{v_j}\wedge \cdots \wedge v_{i+1}).
	\end{align*}
If we wish to emphasize the dependence on $X$, we write $d_{i,X}$ instead of $d_i$.
We also set $d_{-1}=\epsilon$.
By Lemma 3.3 of \cite{Ch16}, the maps $d_i$ are well-defined
$\mbb H$-maps with $d_{i+1}\circ d_i=0$, so \eqref{Koszul-type-complex} is a complex.

	\subsection{Alternate form of the Koszul-type resolution}
	In this section, we provide an alternative form for the differential $d_i$, which involves the terms $\widetilde{v}$ defined in \eqref{Alternate-Kozual-type-resolution} below, which is advantageous for later computations.
	
	For $v \in V$, we define the following element of $\mbb H$:
	\begin{equation}\label{Alternate-Kozual-type-resolution}
		\widetilde{v}=v-\frac{1}{2}\sum_{\alpha \in R^+}k_{\alpha}\langle v, \alpha^{\vee} \rangle t_{s_\alpha}.
	\end{equation}
	
	We consider the maps $\widetilde{d_i}: \mbb H\otimes_{\mbb C[W]}({\rm Res}_WX\otimes \wedge^{i+1} V)\rightarrow \mbb H\otimes_{\mbb C[W]}({\rm Res}_WX \otimes \wedge^iV)$  as follows:
	\begin{equation}
		\begin{split}
			& \widetilde{d_i}(h\otimes (x \otimes v_1\wedge \cdots \wedge v_{i+1}))\\
			=&\sum_{j=1}^{i+1}(h\widetilde{v}_j\otimes x \otimes v_1\wedge \cdots \wedge \widehat{v}_j\wedge \cdots\wedge v_{i+1}-h\otimes \widetilde{v}_j.x\otimes v_1\wedge \cdots\wedge \widehat{v}_j\wedge \cdots \wedge v_{i+1}).
		\end{split}
	\end{equation}	
	This definition coincides with the one given in the previous subsection, so that $\widetilde{d}_i = d_i$.

	We now employ the projective resolution constructed in \eqref{Koszul-type-complex} to build a complex for computing  Ext-groups. Let $X$ and $Y$ be $\mbb H$-modules.
	Applying the contravariant functor ${\rm Hom}_{\mbb H}(\cdot, Y )$ to the projective resolution of $X$, we obtain a sequence of induced differential maps
	\begin{align*}
		d_i^\vee : {\rm Hom}_{\mbb H}(\mbb H \otimes_{\mbb C[W]} ({\rm Res}_W X \otimes \wedge^i V), Y ) \rightarrow {\rm Hom}_{\mbb H}(\mbb H \otimes_{\mbb C[W]} ({\rm Res}_W X \otimes\wedge^{i+1} V ), Y ).
	\end{align*}
	By Frobenius reciprocity, these maps descend to well-defined induced differentials
	\begin{align*}
		d_{i,W}^\vee : {\rm Hom}_{\mbb H}({\rm Res}_W X \otimes \wedge^iV, {\rm Res}_W Y ) \rightarrow {\rm Hom}_\mbb H({\rm Res}_W X \otimes \wedge^{i+1}V, {\rm Res}_W Y).
	\end{align*}
	This procedure yields a complex of homomorphism groups
	\begin{equation}\label{induced-maps}
		\begin{split}
			0 \leftarrow &{\rm Hom}_W ({\rm Res}_W X \otimes \wedge^n V, {\rm Res}_W Y)\overset{d_{n-1,W}^\vee}{\leftarrow} \cdots \overset{d_{i,W}^\vee}{\leftarrow}\\&{\rm Hom}_W({\rm Res}_WX \otimes \wedge^iV, {\rm Res}_W Y)
			\overset{d_{i-1,W}^\vee}{\leftarrow}\cdots \overset{d_{0,W}^\vee}{\leftarrow} {\rm Hom}_W ({\rm Res}_W X, {\rm Res}_W Y ) \leftarrow 0,
		\end{split}
	\end{equation}
The differential $d_{i,W}^\vee$ admits the following explicit formula:
	\begin{equation}\label{zero differ}
		\begin{split}
		d_{i,W}^\vee(\eta)(x \otimes v_1\wedge\cdots \wedge v_{i+1}) =\sum_{j=1}^{i+1}&(-1)^{j+1}\widetilde{v}_j.\eta
		(x\otimes v_1\wedge \cdots\wedge \widehat{v}_j\wedge \cdots \wedge v_{i+1})\\
		-&\sum_{j=1}^{i+1}(-1)^{j+1}
		\eta(\widetilde{v}_j.x \otimes v_1\wedge \cdots \wedge \widehat{v}_j \wedge \cdots \wedge v_{i+1}),
		\end{split}
	\end{equation}
	where the action of $\widetilde{v}_j$ on the term $\eta(x \otimes v_1 \wedge\cdots \wedge \widetilde{v}_j \wedge \cdots \wedge v_{i+1})$ is via the action of $\widetilde{v}_j$ on $X$.
	Thus we obtain the following proposition.
	\begin{proposition}\label{vanish}
		Let $X$ and $Y$ be one $K$-type representations.
        Then ${\rm Ext}_{\mbb H}^i(X,Y)$ is naturally isomorphic
		to the $i$-th homomorphsim group of complex in \eqref{induced-maps}.
	\end{proposition}
	\begin{proof}
		Since $X$ and $Y$ are one $\rm K$-type representations, according to \cite{BM99}, one has $\widetilde{v}_j.x=0$ and $\widetilde{v}_j.\eta(x\otimes v_1\wedge\cdots\wedge\widehat{v}_j\wedge\cdots\wedge v_{i+1})=0$. It follows that $d_{i,W}^\vee=0$, and the proof is completed.
	\end{proof}

By the multiplication of the extension group of ${\rm Ext}_{\mbb H}^j(X,Y)$, we have the following two propositions
\begin{proposition}
 For any two one $K$-type representations $X$ and $Y$ of $\mbb H$, we have
    \begin{align*}
        {\rm Ext}_{\mbb H}^j(X, Y)\cong {\rm Ext}_{\mbb H}^j(Y,X).
    \end{align*}
\end{proposition}	
\begin{proof}	
	Let $\chi_X$ and $\chi_Y$ be the characters of the representations ${\rm Res}_W X$ and ${\rm Res}_W Y$, and $\chi_j$ be the character of the $\wedge^j V$. Then we have
	\begin{align*}
		{\rm dim \ Ext}_{\mbb H}^j(X,Y)=\langle \chi_X \chi_j, \chi_Y \rangle = \frac{1}{|W|}\sum_{w\in W}\chi_X(w)\chi_j(w)\overline{\chi_Y(w)}.
		\end{align*}
		All characters are real-valued, so $\overline{\chi_Y}=\chi_Y$.
	The summand is thus completely symmetric in $\chi_X, \chi_Y$ and $ \chi_j$, which yields
	\begin{align*}
		{\rm dim \ Ext}_{\mbb H}^j(X,Y)={\rm dim \ Ext}_{\mbb H}^j(Y,X).
	\end{align*}
	This establishes the desired isomorphism.
	\end{proof}
	
Let $(\pi, X)$ be a representation of $\mbb H$. Define $(\pi^\ast, X)$ the representation as follows
\begin{align*}
    \pi^\ast (t_w)\cdot x=(-1)^{l(w)} \pi(t_w) \cdot x,  \ \ \pi^\ast(v)\cdot x=-\pi(v)\cdot x, \  \ \forall w\in W, v\in V, x \in X.
\end{align*}

If $X$ is a one $K$-type representation of $\mbb H$, then $X^\ast$ is also a one $K$-type representation, and ${\rm Res}_W X^\ast= {\rm Res}_W X \otimes {\rm sgn}$. We then have the following proposition.
\begin{proposition}
    Let $X$ and $Y$ be two one $K$-type representation of $\mbb H$. Then we have
    \begin{align*}
    {\rm Ext}_{\mbb  H}^{j}(X, Y)={\rm Ext}_{\mbb H}^{n-j}(X^\ast,Y), \ \ \forall j.
    \end{align*}
\end{proposition}

	\section{The Extension group of one $K$-type representations}\label{section 3}

\subsection{The classical group of type $A$}
	
The representation $V$ admits the following decomposition into irreducible constituents:
\begin{align*}
	V \cong L(n) \oplus L(n-1,1).
\end{align*}
Moreover, we have
\begin{align*}
	\wedge^m(L(n-1,1))\cong L(n-m,1^m)
\end{align*}
for each $0\leq m\leq n-1$, where $(n-m, 1^m)$ stands for the $m$-th hook
\begin{align*}
	(n-m, 1^m)=(n-m,\; \overbrace{1,\,1,\,\dots,\,1}^{m \ {\rm times}}).
	\end{align*}
Since $L(n)$ is isomorphic to the trivial representation, it follows that
	\begin{equation}
		\wedge^m V \cong
		\begin{cases}
			L(n) & m=0,\\
			L(n-m+1,\, 1^{\,m-1}) \oplus L(n-m,\, 1^{\,m}) & 0<m<n,\\
			L(1^{\,n}) & m=n.
		\end{cases}
	\end{equation}
	
	\begin{proposition}\label{theorem 1.1.1}
	Assume that $X$ and $Y$ are one $\rm K$-type representations of $\mbb H_{A_{n-1}}$ corresponding to the rectangular partitions $\lambda=(k^d)$ and $\mu=(k'^{d'})$, respectively.
    If $\{k,d\}\neq \{k',d'\}$, then
		\begin{align*}
			{\rm Ext}_{\mbb H_{A_{n-1}}}^j(X, Y)=0, \ \  \text{for all } j.
			\end{align*}
		\end{proposition}
		\begin{proof}
		According to Proposition \ref{vanish} and \cite{Ha15}, the dimension of ${\rm Ext}_{\mbb H}^j(X,Y)$ is equal to cardinality of ${\rm PW}_j(\lambda, \mu)$.
        Assume that $\zeta$	is a partition of $n-j$ such that $\zeta \subset \lambda \cap \mu$ and there exists a picture from $(\lambda/\zeta)^t$ to $\mu/\zeta$.		
        Then we have $(\lambda/\zeta)^t=\mu/\zeta$.

 Suppose that the skew shape $\lambda/\zeta$ has $a$ rows and $b$ columns.
 By transposition, $\mu/\zeta$ then has $b$ rows and $a$ columns.
 Since $\zeta \subset \lambda$, the number $a$ of rows of $\lambda/\zeta$ satisfies $a \le d$, while its number $b$ of columns satisfies $b \le k$.
 Similarly, as $\zeta \subset \mu$, the number $b$ of rows of $\mu/\zeta$ satisfies $b \le d'$, and its number $a$ of columns satisfies $a \le k'$.

Let $s$ be the number of boxes of $\zeta$ in the first row, and $t$ the number in the first column.
Since $\zeta$ is common to both $\lambda$ and $\mu$, we have $s \le \min(k,k')$ and $t \le \min(d,d')$.
The first row of $\lambda/\zeta$ has length $k-s$, and the first column of $\mu/\zeta$ has height $d'-t$. Equality of the two skew shapes implies that the length of the first row of one equals the height of the first column of the other, yielding
\begin{equation}\label{eq:first}
    k-s = d'-t.
\end{equation}

Comparing the first column of $\lambda/\zeta$ with the first row of $\mu/\zeta$, we find that the height of the first column of $\lambda/\zeta$ is $d-t$, while the length of the first row of $\mu/\zeta$ is $k'-s$.
Transposition of the skew shapes gives
\begin{equation}\label{eq:second}
    d-t = k'-s.
\end{equation}

From \eqref{eq:first} and \eqref{eq:second}, we obtain $k-d' = s-t$ and $d-k' = t-s$.
Adding these two identities together:
\begin{align*}
   (k-d')+(d-k')=0 \implies k+d = d'+k'.
\end{align*}
Combined with condition $kd=k'd'=n$, we conclude $\{k,d\}=\{k',d'\}$.
\end{proof}

The following theorem gives a formula for the dimension of Ext between two one $K$-type representations attached to a pair of rectangular partitions.

\begin{theorem}\label{theorem 1.1.2}
Let $X$ and $Y$ be two one $\rm K$-type representations corresponding to the rectangular partitions $\lambda = (k^d)$ and $\mu = (d^k)$, respectively, and set $s={\rm min}(d,k)$. Then
\begin{align*}
{\rm dim \ Ext}_{\mbb H_{A_{n-1}}}^j(X,Y)
=[q^{n-j}]\prod_{i=1}^{s}(1+q^{2i-1}),
\end{align*}
where $[q^m]$ denotes the coefficient of $q^m$.
\end{theorem}
\begin{proof}
	 The intersection $\lambda \cap \mu$ is an $s\times s$ square. Define
	$$\mcal S_m=\{\zeta \vdash m \mid \zeta=\zeta^t, \zeta \subset (s^s)\}$$
	to be the set of all partitions of size $m$
	that are self-conjugate and fit inside an $s\times s$ square.
    For each $\zeta \in \mcal S_{n-j}$, there is exactly one picture in ${\rm PW}(\lambda, \mu; \zeta)$.
    Then it follows that $|{\rm PW}_j(\lambda, \mu)|=|\mcal S_{n-j}|$.
	It is well known that self-conjugate partitions are in bijection with odd partitions via the Frobenius coordinate correspondence.
	The cardinality of $\mcal S_m$ is the coefficient of $q^m$ in the generating function
$$ |\mathcal S_m| = [q^{n-j}]\prod_{i=1}^{s}(1+q^{2i-1}). $$
	In other words, $|S_m|$ counts the number of partitions of $m$ into distinct odd integers chosen from $\{1,3, \dots, 2s-1\}$.
	\end{proof}
	The following proposition gives the dimension of the self-extensions of a one $\rm K$-type representation.
\begin{proposition}
	Let $X$ be a one $\rm K$-type representation corresponding to the rectangular partitions $\lambda=(k^d)$. Then
	\begin{align*}
		{\rm dim \ Ext}_{\mbb H_{A_{n-1}}}^j(X, X)=[q^{j}]\prod_{i=1}^{s}(1+q^{2i-1}).
	\end{align*}
\end{proposition}\label{theorem 1.1.3}
Moreover, we have the following vanishing result for the Euler characteristic of extensions between one $\rm K$-type representations.
\begin{proposition}
Let $X$ and $Y$ be two one $\rm K$-type representations. Then
\begin{align*}
	{\rm EP}(X, Y):=\sum(-1)^j {\rm dim \ Ext}_{\mbb H_{A_{n-1}}}^j(X,Y)=0.
	\end{align*}
	\end{proposition}	
	
	\section{The classical group of type $B$}
	
Elements of $B_n$ are expressed as pairs $(\varepsilon, \sigma)$, where $\varepsilon \in \{\pm1\}^n$ and $\sigma \in S_n$.
The group multiplication is defined by
\begin{align*}
(\varepsilon, \sigma)(\delta, \tau)=(\varepsilon\cdot(\sigma\cdot\delta), \sigma\tau),
\end{align*}
where $(\sigma\cdot \delta)_i=\delta_{\sigma^{-1}(i)}$.
Let $[n]=\{1, \cdots, n\}$.
For positive integers $a, b$ with $a+b=n$, fix the decomposition
\begin{align*}
[a]&=\{1,\dots,a\}, \quad [b]=\{a+1,\dots,n\}.
\end{align*}
Define the subgroups
\begin{align*}
    B_a^b&=B_a\times B_b=\{(\varepsilon, \sigma) \in B_n \mid \sigma([a])=[a], \ \sigma([b])=[b]\}.\\
\end{align*}
The representation $V$ is isomorphic to the irreducible representation $L_{(n-1)(1)}^{B_n}$.
\begin{proposition}
	The exterior power representation $\wedge^mV$ is isomorphic to the irreducible representation $L_{(n-m, 1^m)}^{B_n}$.
	\end{proposition}
\begin{proof}
Let $\{x_{i_1} \wedge \cdots \wedge x_{i_m} \mid i_1< \cdots<i_m\}$ denote the standard basis of $\wedge^m V$.
The representation $L_{(n-m)(1^m)}^{B_n}$ has a basis indexed by $m$-element subsets $C_m \subset [n]$, with elements of the form
\begin{align*}
\prod_{i<j \in C_m}(x_i^2-x_j^2)\prod_{i\in C_m} x_i.
\end{align*}
Define $\Phi$ by sending
\begin{align*}
x_{i_1} \wedge \cdots \wedge x_{i_m}
\longmapsto
\prod_{i<j \in \{i_1,\dots,i_m\}}(x_i^2-x_j^2)\prod_{i\in \{i_1,\cdots, i_m\}}x_i.
\end{align*}
A straightforward computation shows that $\varphi$ intertwines the action of the simple reflections of $B_n$. That is, for every simple reflection $s_i$ and every basis element of $\wedge^m V$, we have
	\begin{align*}
		\Phi\bigl(s_i \cdot (x_{i_1} \wedge \cdots \wedge x_{i_m})\bigr)
		= s_i \cdot \Phi(x_{i_1} \wedge \cdots \wedge x_{i_m}).
	\end{align*}
    Therefore $\Phi$ is an isomorphism.
	\end{proof}	
	
By the Littlewood每Richardson rule in \cite{Fu97}, we have the following lemma.
	
\begin{lemma}\label{multipicity one lemma}
	Let $\nu=(a^k)$ be a rectangle partition and $\lambda \subset \nu$ be a subpartition.
Define \(\nu_\lambda = (\rho_1, \dots, \rho_k)\) with \(\rho_i = a - \lambda_{k+1-i}\).
Then we have $(\nu)_{\lambda}^t=(\nu^t)_{\lambda^t}$ and
$C_{\lambda,\mu}^{\nu} = \delta_{\mu, \nu_\lambda}$.
	\end{lemma}
	
The sign representation of $B_n$ is isomorphic to $L_{(\emptyset)(1^n)}^{B_n}$.
\begin{proposition}
For any irreducible representation $L_{(\lambda)( \mu)}^{B_n}$, we have
		\begin{align*}
			L_{(\lambda)( \mu)}^{B_n} \otimes L_{(\emptyset)(1^n)}^{B_n}=L_{(\mu^t)( \lambda^t)}^{B_n}.
		\end{align*}
	\end{proposition}
	\begin{proof}
Let $\lambda \vdash a$ and $\mu \vdash b$ with $a+b=n$.
Recall that
$$L_{(\lambda)(\mu)}^{B_n} \cong \operatorname{Ind}_{B_a \times B_b}^{B_n} (\rho_\lambda^0 \boxtimes \rho_\mu^1),$$
where $\rho_\lambda^0$ and $\rho_\mu^1$ denote representations with trivial and sign action on the
$(\mbb Z_2)^n$-factor, respectively.
The sign representation \(L_{\emptyset,(1^n)}^{B_n}\) restricts as
\begin{align*}
L_{(\emptyset)(1^n)}^{B_n}|_{B_a \times B_b}
= L_{(\emptyset)(1^a)}^{B_a} \boxtimes L_{(\emptyset)(1^b)}^{B_b}.
\end{align*}
By the projection formula,
\begin{align*}
L_{(\lambda)(\mu)}^{B_n} \otimes L_{(\emptyset)(1^n)}^{B_n}
&= {\rm Ind}_{B_a \times B_b}^{B_n}
\Bigl( (\rho_\lambda^0 \boxtimes \rho_\mu^1)
\otimes (L_{(\emptyset)(1^a)}^{B_a} \boxtimes L_{(\emptyset)(1^b)}^{B_b}) \Bigr) \\
&={\rm Ind}_{B_a \times B_b}^{B_n}
\Bigl((\rho_\lambda^0 \otimes L_{(\emptyset)(1^a)}^{B_a})
\boxtimes
(\rho_\mu^1 \otimes L_{(\emptyset)(1^b)}^{B_b}) \Bigr) \\
&=\operatorname{Ind}_{B_a \times B_b}^{B_n}
(\rho_{\lambda^t}^1 \boxtimes \rho_{\mu^t}^0)
\cong L_{(\mu^t)(\lambda^t)}^{B_n}.
\end{align*}
The penultimate equality uses the standard fact that tensoring with the sign representation transposes the partition.
\end{proof}

The double coset space $B_a^b\backslash G/B_k^l$ appears naturally in the Mackey decomposition for restrictions of representations of $B_n$ to $B_k^l$.

\begin{proposition}\label{mackey decomposition}
 The double cosets $B_a^b \backslash B_n / B_k^l$ are classified by quadruples $(a_1,a_2,b_1,b_2) \in \mathbb Z_{\ge 0}^4$ satisfying
\begin{equation}\label{addition}
a_1 + b_1 = k,\quad a_2 + b_2 = l,\quad a_1 + a_2 = a,\quad b_1 + b_2 = b.
\end{equation}
For a double coset represented by $s$ and classified by $(a_1,a_2,b_1,b_2)$,  we have
\begin{align*}
B_a^b \cap sB_k^l s^{-1}\cong B_{a_1}\times B_{b_1}\times B_{a_2} \times B_{b_2}.
\end{align*}
\end{proposition}
\begin{proof}
We first identify the double coset space. Define
\begin{align*}
\Phi : &B_a^b\backslash B_n/B_k^l \;\longleftrightarrow\; (S_a \times S_b)\backslash S_n/(S_k \times S_l), \\
&\Phi (B_a^b(\varepsilon, \sigma)B_k^l):=(S_a\times S_b)\sigma (S_k\times S_l).
\end{align*}
This is well-defined. Indeed, if $B_a^b(\varepsilon, \sigma)B_k^l=B_a^b(\varepsilon', \sigma')B_k^l$, then
$(\varepsilon',\sigma')=h(\varepsilon, \sigma)k$ with $h \in B_a^b, k\in B_k^l$.
Taking permutation parts gives $\sigma'=\tau\sigma\rho$ with $\tau \in S_a\times S_b, \rho\in S_k\times S_l$, so $\sigma$ and $\sigma'$ lies in the same double coset.

Take any $\sigma \in S_n$ and set $g=(1, \sigma) \in B_n$. Then $\Phi(B_a^bgB_k^l)$ is the given double coset, so surjectivity follows.

For injectivity, suppose $\Phi(B_a^b(\varepsilon,\sigma)B_k^l)=\Phi(B_a^b(\varepsilon',\sigma')B_k^l)$.
Then we have $\sigma'=\tau\sigma\rho$ for some $\tau\in S_a\times S_b, \rho\in S_k\times S_l$.
Choose lifts $h=(\varepsilon'\cdot(\tau\cdot\varepsilon)^{-1}\cdot(\tau\sigma\cdot\theta)^{-1},\tau)\in B_a^b$ and $k=(\theta,\rho)\in B_k^l$. Then \(h(\varepsilon,\sigma)k=(\varepsilon',\sigma')\), proving equality of the \(B_n\)-double cosets. Hence \(\Phi\) is bijective.

It remains to classify the double cosets $(S_a\times S_b)\backslash S_n/(S_k\times S_l)$. For a representative $\sigma\in S_n$, define
\begin{equation*}
a_1=|[a]\cap\sigma[k]|,\quad a_2=|[a]\cap\sigma[l]|,\quad
b_1=|[b]\cap\sigma[k]|,\quad b_2=|[b]\cap\sigma[l]|.
\end{equation*}
The equations \eqref{addition} follow immediately.
If $\sigma'=\tau\sigma\rho$ with $\tau\in S_a\times S_b, \rho\in S_k\times S_l$, then the four intersection cardinalities are invariant, so the map is well-defined.

 Suppose two double cosets give the same quadruple. Let
 $$A=\sigma([k]), \ B=\sigma([l]), \  A'=\sigma'([k]), \ B'=\sigma'([l]).$$
 Then there exists $\tau \in S_k\times S_l$ mapping mapping each of $A\cap [a], A\cap [b]$
to the corresponding intersections of $A'$, and similarly for $B$.
Defining $\rho=\sigma^{-1}\tau^{-1}\sigma'$, we obtain $\tau\sigma\rho=\sigma'$.
Since $\rho([a])=[a]$ and $\rho([b])=[b]$, it follows that $\rho\in S_a\times S_b$.
Hence the two double cosets coincide.
Surjectivity follows by choosing a subset $A\subset[n]$ of size $a$ with $|A\cap[k]|=a_1$ and $|A\cap[l]|=a_2$, then taking any bijection $\sigma$ with $\sigma([a])=A$.

It remains to compute the intersection $B_a^b \cap s^{-1}B_k^ls$.
For $s=(\varepsilon,\sigma)\in B_n$, an element $g=(\eta,\tau)$ lies in this intersection if and only if
$\tau$ preserves each of the sets $[a], [b], s[k]$ and $s[l]$. Define
$$A_1=[a]\cap s[k], \ A_2=[a]\cap s[l], \ B_1=[b]\cap s[k], \ B_2=[b]\cap s[l].$$
These four sets form a partition of $[n]$ with cardinalities $a_1, a_2, b_1, b_2$, respectively.
Since $\tau$ preserves each of the original four sets, it also preserves each of $A_1, A_2, B_1, B_2$.
Then
\begin{align*}
B_a^b\cap s B_k^l s^{-1}\cong B_{a_1}\times B_{b_1}\times B_{a_2}\times B_{b_2}.
\end{align*}
This completes the proof.
\end{proof}

\subsection{The one partition}
	
	In this subsection, we consider the case where $X$ is a one $\rm K$-type representation corresponding to a bipartition $(\sigma_L)(\emptyset)$.

	\begin{proposition}\label{both one partition}
	Let $Y$ be a one $\rm K$-type representation of $\mbb H_{B_n}$, and the corresponding bipartition is $(\lambda)(\emptyset)$.
	Then we have
	\begin{align*}
		{\rm dim \ Ext}_{\mbb H_{B_n}}^j (X, Y)=
		\begin{cases}
			1, &  \text{if} \  \sigma_L=\lambda, j=0.\\
			0, &   \text{otherwise}.
		\end{cases}
		\end{align*}
		\end{proposition}
	\begin{proof}
Since $Y|_{\mbb C[W]}$ is $L_{(\lambda)(\emptyset)}^{B_n}$ or $L_{(\emptyset)(\mu)}^{B_n}$, we have
\begin{equation}\label{multiplicity one}
\begin{split}
L_{(\sigma_L)(\emptyset)}^{B_n} \otimes L_{(n-m)(1^m)}^{B_n}&=\rho_{\sigma_L}^0 \otimes {\rm Ind}_{B_{n-m}\times B_m}^{B_n} \rho_{(n-m)}^0 \boxtimes \rho_{(1^m)}^1\\
&=\bigoplus_{|\xi|=m} c_{\zeta, \xi}^{\sigma_L} {\rm Ind}_{B_{n-m}\times B_m}^{B_n} \rho_{\zeta}^0\boxtimes \rho_{\xi^t}^1=\bigoplus_{|\xi|=m} c_{\zeta, \xi}^{\sigma_L} L_{(\zeta)(\xi^t)}^{B_n}.
\end{split}
\end{equation}
Then ${\rm dim \ Ext}_{\mbb H_{B_n}}^j(X,Y)=1$ if $\lambda=\sigma_L, j=0$ or if $\mu=\sigma_L, j=n$.
 In all other cases, the Ext groups vanish.
\end{proof}
Let $\alpha=(a^p)(b^q)$ be a bipartition of $n$ with $p-a=q-b+\delta$.  Set
$$x=\min\{a,b\}, \ \  y=\min\{p,q\}, \ \ z=a+b-\delta_{p,y}b-\delta_{q,y}a.$$
Denote by $S_{\alpha}$ the set of partitions of $n$ such that
\begin{align*}
	S_{\alpha}=\Big\{ a+b-\nu_y, \cdots a+b-\nu_1, z^{p+q-2y}, \nu_1,\cdots, \nu_1 \ \Big| \ \nu \subset (x^y) \Big\}.
	\end{align*}
	
\begin{proposition} \label{one-two bipartition}
 Let $Y$ be a one $\rm K$-type representation of $\mbb H_{B_n}$, and the  corresponding bipartition is $(a^p)(b^q)$.
 Then we have
 \begin{align*}
		{\rm dim \ Ext}_{\mbb H_{B_n}}^j (X, Y)=
		\begin{cases}
			1, &  \text{if} \  \sigma_L\in S_{(a^p)(b^q)} \ \text{and} \ j=bq,\\
			0, &   \text{otherwise}.
		\end{cases}
		\end{align*}
\end{proposition}
\begin{proof}
By the symmetry of the Littlewood--Richardson coefficients $c_{\alpha, \beta}^{\nu}=c_{\beta, \alpha}^{\nu}$,  we may assume without loss of generality that $p\geq q$.
From the arguments in \cite{Ok98}, we have
\begin{align*}
    {\rm Ind}_{B_{ap}\times B_{bq}}^{B_n} L_{(a^p)(\emptyset)}^{B_{ap}}\boxtimes L_{(b^q(\emptyset)}^{B_{bq}}=\bigoplus_{\lambda\in S_\alpha} L_{(\lambda)(\emptyset)}^{B_n}.
\end{align*}
The result then follows from equation \eqref{multiplicity one}.
\end{proof}

\subsection{The bipartition cases}
Let $\alpha=(a^p)(b^q)$ and $\beta=(c^k)(d^l)$ be two rectangular bipartitions of $n$.
We define the following set
\begin{equation}\label{bipartition complement}
	\begin{split}
    \Pi_{\alpha,\beta}^{j}=\{(\zeta, \xi) \mid \zeta \subset (a^p)\cap(l^d), \  \xi \subset (b^q)\cap (k^c),  \ (a^p)_{\zeta}&=(c^k)_{\eta},  \\  (q^b)_{\eta}&=(l^d)_{\zeta}, \ |\zeta|+|\eta|=j\}.
\end{split}
\end{equation}

\begin{theorem}\label{bipartition multiplication of B}
    Let $X$ and $Y$ be two one $\rm K$-type representations of $\mbb H_{B_n}$ with associated bipartitions $\alpha=(a^p)(b^q)$ and $\beta=(c^k)(d^l)$, respectively. Then ${\rm dim\ Ext}_{\mbb H_{B_n}}^j(X,Y)=0$ for all $j$ if and only if none of the following conditions holds:
    \begin{align}\label{four conditions}
    	\begin{array}{cl}
    		\mathrm{I}   & a=c,\ p=k,\ b=d,\ q=l\\[2pt]
    		\mathrm{II}  & a=l,\ p=d,\ c=q,\ k=b\\[2pt]
    		\mathrm{III} & a=c=q=l,\quad p+b=k+d=2a+\tau\\[2pt]
    		\mathrm{IV}  & p=b=k=d,\quad a+q=c+l=2p-\tau
    	\end{array}
    	 \end{align}
  	If at least one of the conditions holds, set
  \begin{align*}
  \begin{aligned}
  	R&=\min(a,c),&
  	S&=\min(p,k),\\
  	W&=\min(a,c,q,l),&
  	H&=\min(p,k,b,d),
  \end{aligned}
  \end{align*}
  	and
  \begin{equation*}
  	j_0=ap+ck-2RS.
  \end{equation*}
  Then
  \begin{equation*}
  		 {\rm dim \ Ext}_{\mbb H_{B_n}}^j(X,Y)=
  			[q^j] q^{j_0} [W+H, H]_{q^2}, \quad  \text{where} \ [W+H, H]_{u} =\prod_{i=1}^H \frac{1-u^{W+i}}{1-u^i}.
  \end{equation*}
\end{theorem}
\begin{proof}
According to Proposition \ref{mackey decomposition} and Mackey theorem, we have
\begin{align*}
  L_{(\gamma)(\delta)}^{B_n}|_{B_{n-j} \times B_{j}}=\bigoplus_{\substack{  |\xi|+|\tau|=j\\ |\zeta|+|\nu|=n-j}}c_{\zeta, \xi}^{\gamma}c_{\nu, \tau}^{\delta} L_{(\zeta)(\nu)}^{B_{n-j}}\boxtimes L_{(\xi)(\tau)}^{B_{j}}.
\end{align*}
Then
\begin{equation}\label{tensor product multiplication}
\begin{split}
  L_{(\gamma)(\delta)}^{B_n} \otimes \wedge^jV&=\bigoplus_{\substack{  |\xi|+|\tau|=j\\ |\zeta|+|\nu|=n-j}}c_{\zeta, \xi}^{\gamma}c_{\nu, \tau}^{\delta}{\rm Ind}_{B_{n-j}\times B_j}^{B_n} (L_{(\zeta)(\nu)}^{B_{n-j}}\boxtimes L_{(\xi)(\tau)}^{B_{j}})\otimes (\rho_{(n-j)}^0 \boxtimes \rho_{(1^j)}^1)\\
  &=\bigoplus_{\substack{  |\xi|+|\tau|=j\\ |\zeta|+|\nu|=n-j}}c_{\zeta, \xi}^{\gamma}c_{\nu, \tau}^{\delta} {\rm Ind}_{B_{n-j}\times B_j}^{B_n} L_{(\zeta)(\nu)}^{B_{n-j}}\boxtimes L_{(\tau^t)(\xi^t)}^{B_{j}}\\&=\bigoplus_{\substack{  |\xi|+|\tau|=j\\ |\zeta|+|\nu|=n-j}}c_{\zeta,\xi}^{\gamma}c_{\nu,\tau}^{\delta}c_{\zeta, \tau^t}^{\lambda}c_{\nu, \xi^t}^{\mu} L_{(\lambda)(\mu)}^{B_n}.
  \end{split}
\end{equation}
Since the differentials in \eqref{zero differ} are zero, we have
\begin{equation}\label{Littlewood每Richardson coefficients}
 {\rm dim \ Ext}_{\mbb H_{B_n}}^j(X,Y)=\sum_{\substack{  |\xi|+|\tau|=j\\ |\zeta|+|\nu|=n-j}}c_{\zeta,\xi}^{(a^p)}c_{\nu,\tau}^{(b^q)}c_{\zeta, \tau^t}^{(c^k)}c_{\nu,\xi^t}^{(d^l)}.
\end{equation}
From the Lemma \ref{multipicity one lemma}, all Littlewood-Richardson coefficients appearing in \eqref{Littlewood每Richardson coefficients} are at most $1$, and the sum coincides with the cardinality of $\Pi_{\alpha,\beta}^{j}$.

We now prove that there exist $j$ such that $\Pi_{\alpha,\beta}^{j}$ is not empty if and only one of the conditions in \eqref{four conditions} holds.
Assume that a pair $(\zeta,\eta)$ satisfying \eqref{bipartition complement} exists. Put
\begin{align}\label{conditionf}
h=k-p,\qquad g=b-d, \qquad s=c-a,\qquad t=q-l.
\end{align}
By the Proposition \ref{conditionBCD}, we have
$(k-p)-(b-d)=(c-a)-(q-l)$, that is,
\begin{align}\label{conditiona}
	h-g=s-t.
	\end{align}
Denote this common value by $\gamma$.  We will show that $\gamma=0$.
From
\begin{align*}
(a^p)_\zeta=(c^k)_\eta,\qquad \text{for} \  \max(1,h+1) \leq i \leq k,
\end{align*}
we compare the $(k+1-i)$-th rows of the two sides and obtain
\begin{equation}\label{conditionb}
	\eta_i=s+\zeta_{i-h}.
\end{equation}
The requirement $i\ge h+1$ guarantees $i-h\ge1$, so that no illegal row index occurs.
Similarly, from
\begin{align*}
(l^d)_\zeta=(q^b)_\eta,\qquad \text{for} \ \max(1,g+1) \leq i \leq b,
\end{align*}
we get
\begin{equation}\label{conditionc}
	\eta_i=t+\zeta_{i-g}.
\end{equation}
Assume $\tau>0$. Then $h>g$ and $s>t$.
We distinguish two cases.

\noindent\emph{Case A: $h<b$.}
Take $i=\max(1,h+1)$.
This row index is admissible in both \eqref{conditionb} and \eqref{conditionc}.
Applying them simultaneously gives
$s+\zeta_{i-h}=t+\zeta_{i-g}$.
Rearranging the terms, we obtain
\begin{equation}\label{conditiond}
	s-t=\zeta_{i-g}-\zeta_{i-h}.
\end{equation}
On the other hand, $h>g$ implies
$i-h<i-g$, and since the rows of a partition are weakly decreasing, it follows that
$\zeta_{i-h}\ge\zeta_{i-g}$.
Combining this with \eqref{conditiond} yields
$s-t\le0$,  which contradicts $s-t=\delta>0$.

\noindent\emph{Case B: $h\ge b$.}
Since $\eta\subseteq(q^b)$, the partition $\eta$ has at most $b$ rows, and hence $\eta_{h+1}=0$.
Comparing the $p$-th rows in $(a^p)_\zeta=(c^k)_\eta$,
The assumption $h\ge b>0$ gives $k=p+h>p$. Comparing the $p$-th rows in $(a^p)_\zeta=(c^k)_\eta$, we get
$$a-\zeta_1=c-\eta_{k+1-p}=c-\eta_{h+1}=c.$$
Thus
\begin{equation*}
	\zeta_1=a-c.
\end{equation*}
For the second complement identity, the $d$-th row of the left-hand side is
\begin{equation*}
		((l^d)_\zeta)_d =l-\zeta_1=l-a+c=q+(c-a-q+l)=q+\tau>q.
\end{equation*}
However, the right-hand side $(q^b)_\eta$ is contained in the rectangle $(q^b)$, so each of its rows is at most $q$, contradicting the equality of the two partitions.
Hence $\tau>0$ is impossible, whether $h<b$ or $h\ge b$.

 If $\tau<0$, interchanging the two complement identities in \eqref{bipartition complement} swaps $h$ with $g$ and $s$ with $t$,
Then one cannot have $h>g$ and $s>t$ simultaneously.
 Hence $\delta<0$ is also impossible.
Therefore $\tau=0$. Together with \eqref{conditiona} we obtain
\begin{equation}\label{conditione}
		k-p=b-d,\qquad c-a=q-l.
\end{equation}

We now prove that non-emptiness forces one of the four conditions to hold. If $h\neq 0$, then $c=q$ and $a=l$.
Likewise, if $s\ne0$, then $p=d$ and $k=b$.
We give a classification according to whether $h$ and $s$ vanish.
Only four possibilities remain.

\noindent\textbf{Case I: $h=0$, $s=0$.}
By \eqref{conditionf}, we have
$$\boxed{\mathrm{I}:\quad a=c,\ p=k,\ b=d,\ q=l.}$$

\noindent\textbf{Case II: $h\ne0$, $s\ne0$.}
Based on the above analysis, we obtain
$$\boxed{\mathrm{II}:\quad a=l,\ p=d,\ c=q,\ k=b.}$$

\noindent\textbf{Case III: $h\ne0$, $s=0$.}
Here $s=0$ gives $c=a$ and $q=l$, while $h\ne0$ gives $a=l$ and $c=q$. Hence $a=c=q=l$. Substituting $q=a$ and $c=l=a$ into relations
$p-a=q-b+\delta$ and $k-c=l-d+\delta$, respectively, yields
 $p+b=k=d=2a+\tau$. Therefore
 \begin{align*}
 \boxed{\mathrm{III}:\quad a=c=q=l,\qquad p+b=k+d=2a+\delta.}
 \end{align*}

\noindent\textbf{Case IV: $h=0$, $s\ne0$.}
From $h=0$ we have $k=p$ and $b=d$, while $s\ne0$ gives $p=d$
and $k=b$. Hence $p=b=k=d$. Similarly, we obtain
$a+q=c+l=2p-\delta$. Thus
 \begin{align*}
\boxed{\mathrm{IV}:\quad p=b=k=d,\qquad a+q=c+l=2p-\delta.}
 \end{align*}
This completes the necessity.
\begin{equation}
		\exists j\ge0:\ N_j\ne0
		\quad\Longrightarrow\quad
		\text{at least one of }\mathrm{I},\mathrm{II},\mathrm{III},\mathrm{IV}
		\text{ holds}.
\end{equation}

We construct, case by case, a pair $(\zeta,\eta)$ satisfying \eqref{bipartition complement}.  We use the convention
$$(a^0)=\emptyset,\qquad (0^p)=\emptyset.$$

\noindent\textbf{Case I.} Assume $a=c, p=k, b=d$ and $q=l$.
Take $\zeta=\eta=\varnothing$.
Both rectangle containment conditions are clearly satisfied, and
\[
(a^p)_{\emptyset}=(a^p)=(c^k)=(c^k)_{\emptyset}, \qquad (l^d)_{\emptyset}=(l^d)=(q^b)=(q^b)_{\emptyset}.
\]
Thus we obtain a solution of total size $j=0$, and $\Pi_{\alpha,\beta}^{0}\neq \emptyset$.

\noindent\textbf{Case II.} Assume $a=l, p=d, c=q$ and $k=b$.
Take the full rectangles $\zeta=(a^p)$ and $\eta=(c^k)$.
Since
$$(a^p)=(l^d),\qquad (c^k)=(q^b),$$
the containment conditions hold.
 The complement of a full rectangle inside itself is empty, so
\begin{align*}
(a^p)_\zeta=\emptyset=(c^k)_\eta, \qquad (l^d)_\zeta=\emptyset=(q^b)_\eta.
\end{align*}
We obtain a solution of total size $j=ap+ck$,  hence $\Pi_{\alpha,\beta}^{n}\neq \emptyset$.

\noindent\textbf{Case III.}
Assume $a=c=q=l$ and $p+b=k+d=2a+\delta$.
In particular,
\begin{align}\label{conditionh}
	p-k=d-b.
	\end{align}
 Let $r=\min(p,k)$ and define
\begin{equation}\label{conditiong}
		\zeta=(a^{p-r}),\qquad \eta=(a^{k-r}).
\end{equation}
 Clearly $\zeta$ has width $a$ and at most $p$
rows, and $\eta$ has width $a$ and at most $k$ rows.
It remains to check
\begin{align*}
p-r\le d,\qquad k-r\le b.
\end{align*}
By \eqref{conditionh},
\begin{align*}
p-r=\max(p-k,0)=\max(d-b,0)\le d,
\end{align*}
and
\begin{align*}
k-r=\max(k-p,0)=\max(b-d,0)\le b.
\end{align*}
Therefore $\zeta\subseteq(a^p)\cap(l^d)$ and $\eta\subseteq(c^k)\cap(q^b)$.
A partition with $p-r$ full rows has, in
the rectangle $(a^p)$, a complement of $r$ full rows. Hence
$(a^p)_\zeta=(a^r)=(a^k)_\eta$,  so the first identity holds.

Since $\zeta=(a^{p-r})$ has exactly $p-r$ rows,
each of length $a$, the row formula gives
$(a^d)_\zeta=(A^{d-p+r})$ and $(a^b)_\eta=(a^{b-k+r})$.
These row counts are nonnegative by the containment check above.  Since $p+b=k+d$, $(a^d)_\zeta=(a^b)_\eta$.
Thus \eqref{conditiong} yields a solution of total size $j=a|p-r|$.
Hence $\Pi_{\alpha,\beta}^{j}\neq \emptyset$.

\noindent\textbf{Case IV.} Assume $p=b=k=d$
and $a+q=c+l=2p-\delta$. This yields
\begin{equation}\label{conditioni}
	a-c=l-q.
\end{equation}
Let $u=\min(a,c)$ and define
\begin{equation}\label{conditionj}
		\zeta=((a-u)^p),\qquad
		\eta=((c-u)^p).
\end{equation}
Both partitions have at most $p$ rows.
For the widths, we have
\begin{align*}
a-u=\max(a-c,0)=\max(l-q,0)\le l,
\end{align*}
 hence $\zeta\subset (a^p)\cap(l^p)$.
Similarly, we have $c-u=\max(c-a,0)=\max(q-l,0)\le q$.
  Hence $\eta\subseteq(c^p)\cap(q^p)$.
 Computing row by row, we have $(a^p)_\zeta=u^{p}=(c^p)_\eta$.

On the other hand, it holds that
\begin{align*}
(l^p)_\zeta=((l-a+u)^p), \ \ (q^p)_\eta=((q-c+u)^p)
\end{align*}
It follows from \eqref{conditioni} that the two complements are equal.
Thus \eqref{conditionj} yields a solution of total size $j=p|a-c|$.

Let $\mu=(a^p)_\zeta=(c^k)_\eta$. Then
\begin{align*}
\mu\subset(R^S),
\qquad R=\min(a,c),\quad S=\min(p,k).
\end{align*}
Define
\begin{align*}
\nu=(R^S)_\mu.
\end{align*}
Conversely, given $\nu$, one recovers uniquely
	\begin{align*}
		\mu=(R^S)_\nu, \qquad
		\zeta=(a^p)_\mu, \qquad
		\eta=(c^k)_\mu.
	\end{align*}
Under the there parameter types, all the original constraints are equivalent to
\begin{align*}
\nu\subseteq(W^H).
\end{align*}
Specifically, in the first case one has $\zeta=\eta=\nu$, and the remaining constraints give $W,H$ directly.
In the second case the two complement
identities coincide, so every $\nu\subseteq(R^S)$ works.
In the third case all widths are equal,  and using $p+b=k+d$ the
restriction on the number of rows reduces exactly to $\ell(\nu)\le H$.
In the fourth case all heights are equal, and using $a+q=c+l$ the restriction on the widths reduces exactly to $\nu_1\le W$.

Thus we obtain the bijection
\begin{align*}
	\bigcup_{j\ge0}\Pi_{\alpha,\beta}^{j}
	\ \longleftrightarrow\
	\{\nu:\nu\subseteq(W^H)\}.
\end{align*}
The sizes satisfy
\begin{align*}
|\zeta|+|\xi|=ap+ck-2|\mu|=ap+ck-2(RS-|\nu|)=j_0+2|\nu|.
\end{align*}
Therefore
\begin{align*}
	\sharp \Pi_{\alpha,\beta}^{j}=
	\#\Bigl\{
	\nu\vdash\tfrac{j-j_0}{2}:
	\nu_1\le W,\ \ell(\nu)\le H
	\Bigr\},
\end{align*}
with the convention that partitions of non-integral or negative size are counted as zero.

Finally, applying the standard generating function for partitions inside a rectangle,
\begin{align*}
\sum_{\nu\subseteq(W^H)}u^{|\nu|}
=[W+H,H]_u,
\end{align*}
 we complete the proof of the theorem,
 This standard identity and its product expression can be found in  \cite{St12}.
\end{proof}
\section{The classical group of type $D$}

Let $L_{(\lambda)(\mu)}^{B_n}$ be the irreducible representation of $B_n$ indexed by the bipartition $(\lambda)(\mu)$.
Define a representation $L_{[\lambda][\mu]}^{D_n}$ of $D_n$ by restriction:
\begin{align*}
    L_{[\lambda][\mu]}^{D_n} :={\rm Res}_{D_n}^{B_n} L_{(\lambda)(\mu)}^{B_n}.
\end{align*}

The following gives a complete classification of irreducible representations.
 \begin{itemize}
\item $(a)$ The representation $L_{[\lambda][\mu]}^{D_n}$ is irreducible if and only $\lambda\neq \mu$.
\item $(b)$ If $\lambda =\mu$, then $L_{[\lambda][\lambda]}^{D_n}$ decomposes as a direct sum of two inequivalent irreducible representations, denoted $L_{\lambda^+}^{D_n}$ and $L_{\lambda^-}^{D_n}$.
\item $(c)$ One has $L_{[\lambda][\mu]}^{D_n}=L_{[\mu][\lambda]}^{D_n}$.
\end{itemize}	
The representations obtained in $(a)$ and $(b)$ form a complete set of pairwise inequivalent irreducible representations of $D_n$.

Similarly, restricting the wedge module to $D_n$, we have $\wedge^mV \cong L_{[n-m][1^m]}^{D_n}$.

\subsection{The one partition case}

In this subsection, we restrict our attention to extension groups ${\rm Ext}_{\mbb H_{D_n}}^j(X,Y)$ where $X$ corresponds to the bipartition $[\upsilon_L][\emptyset]$.

\begin{proposition}
    Let $Y$ be a one $\rm K$-type representation of $\mbb H_{D_n}$ labeled by the bipartition $[\lambda][\emptyset]$.
    Then
    \begin{align*}
        {\rm dim \ Ext}_{\mbb H_{D_n}}^j(X,Y)=\begin{cases}
			1, &  \text{if} \  \upsilon_L=\lambda, j=0 \ \text{or} \ \upsilon_L=\lambda^t, j=n,\\
			0, &   \text{otherwise}.
		\end{cases}
    \end{align*}
\end{proposition}
\begin{proof}
    For any two representations $V, W$ of $B_n$, we have
    \begin{align}\label{Res multiplication}
        {\rm Res}_{D_n}^{B_n}(V\otimes W)={\rm Res}_{D_n}^{B_n}(V) \otimes {\rm Res}_{D_n}^{B_n}(W).
    \end{align}
By the \eqref{tensor product multiplication},  we have
\begin{align}\label{D_n multiplication}
    L_{[\lambda][\emptyset]}^{D_n}\otimes \wedge^jV=(\rho_{\lambda}^0\otimes \wedge^jV)|_{D_n}=\bigoplus_{|\xi|=j}c_{\zeta, \xi^t}^{\lambda}L_{[\zeta][\xi]}^{D_n}.
\end{align}
By the symmetry $L_{[\zeta][\xi]}^{D_n}=L_{[\xi][\zeta]}^{D_n}$, we have ${\rm dim \ Ext}_{\mbb H}^j(X,Y)=1$, if  $\upsilon_L=\lambda, j=0$ or $\upsilon_L=\lambda^t, j=n$, and vanishes otherwise.
\end{proof}
Using the tensor product decomposition of irreducible $D_n$-representations in \eqref{D_n multiplication}, we derive the following corollary.
\begin{corollary}
    Let $Y$ be a one $\rm K$-type representation of $\mbb H_{D_n}$ such that $Y|_{\mbb C[D_n]}=L_{[a^p][b^q]}^{D_n}$ or $Y|_{\mbb C[D_n]}=L_{[c^k]^{\pm}}^{D_n}$.
    Then
    \begin{align*}
     {\rm dim \ Ext}_{\mbb H}^j(X, Y)=  \begin{cases}
			1, &  \text{if} \ Y|_{\mbb C[D_n]}=L_{[a^p][b^q]}^{D_n},  \lambda \in S_{(a^p)(b^q)} \  \text{and} \ j=ap \ \text{or} \  j=bq,\\
            1, &  \text{if} \ Y|_{\mbb C[D_n]}=L_{[c^k]^{\pm}}^{D_n},  \lambda \in S_{(c^k)(c^k)} \  \text{and} \ j=n/2,\\
			0, &   \text{otherwise}.
		\end{cases}
    \end{align*}
\end{corollary}
\subsection{The two bipartition case}
In this subsection, we consider the case where the single K-type representation $X$ satisfies $X|_{\mathbb C[D_n]} = L_{[a^p][b^q]}^{D_n}$, with $(a^p) \neq (b^q)$.

By Proposition \ref{bipartition multiplication of B} and equation \eqref{Res multiplication}, we obtain the following proposition.

\begin{proposition}\label{theorem 1,3.1}
    Let $X$ and $Y$ be one $K$-type representations of $\mathbb H_{D_n}$ with corresponding bipartitions $\alpha = (a^p)(b^q)$ and $\beta = (c^k)(d^l)$, respectively.
    Then
    \begin{align*}
        \dim \operatorname{Ext}_{\mathbb H}^j(X,Y) = \sharp\,\Pi_{\alpha, \beta}^{j} + \sharp\,\Pi_{\alpha, \overline{\beta}}^{j},
    \end{align*}
    where $\overline{\beta}$ denotes the bipartition $(d^l)(c^k)$.
\end{proposition}

\subsection{The split irreducible representation case}
In this subsection, we consider the  case that the one $\rm K$-type representation $X$ satisfies $X|_{\mathbb C[D_n]} = L_{[a^p]^+}^{D_n}$.
We set $m = ap$ and $n = 2m$.

Denote by $\chi_{\lambda}$ and $\chi_{\lambda^\pm}$  the characters of the representations $L_{[\lambda][\lambda]}^{D_n}$ and $L_{[\lambda]^\pm}^{D_n}$, respectively. The difference character $\Delta_{\lambda}$ is defined by
\begin{align*}
   \Delta_{\lambda}=\chi_{\lambda^+}- \chi_{\lambda^-}.
\end{align*}
For the irreducible representation $\wedge^j V$, denote its character by $\chi_j$.
Fix a element $s \in B_{n}\backslash D_{n}$, we have
\begin{align*}
    \Delta_{\lambda}^s(g):=\Delta_\lambda(sgs^{-1})=-\Delta_\lambda(g), \  \ \chi_j^s(g):=\chi_j(sgs^{-1})=\chi_j(g),  \ \ \forall g \in D_n.
\end{align*}

Let $\tau=(\tau_1,\tau_2,\cdots,\tau_l)$ be a partition. For any positive integer $j$, define
\begin{align*}
    m_j(\tau):=\sharp\{i \mid \tau_i=j \ \}, \ \ z^\tau=\prod_{j}j^{m_j(\tau)} m_j(\tau)!.
\end{align*}

\begin{proposition}\label{theorem 1.3.2}
    If $Y$ is a one $\rm K$-type representation of $\mathbb H_{D_n}$ such that $Y|_{\mathbb C[D_n]} = L_{[c^k]^\theta}^{D_n}$, then
     \begin{align*}
     {\rm dim \ Ext}_{\mbb H_{D_n}}^j(X, Y)=  \begin{cases}
			\frac{1}{2}([q^r][2s,s]_q+\theta(-1)^rN_{m-r}), &  \text{if} \ j=2r, \ (a^p)=(c^k), \\
            \frac{1}{2}([q^{m-r}][2s,s]_q+\theta(-1)^rN_{m-r}), &  \text{if} \ j=2r, \ (a^p)=(k^c), \\
			0, &   \text{otherwise},
		\end{cases}
    \end{align*}
where $s=\min \{a,p\}$ and $N_r=[q^r]\prod_{i=1}^s (1+q^{2i-1})$.
\end{proposition}
\begin{proof}
By the arguments in \cite{Ta15}, we have the decomposition of the tensor product.
\begin{align*}
   \Delta_\lambda \otimes \chi_j
    =T_j(\lambda,\mu) \Delta_\mu, \ \ T_j(\lambda,\mu)=\sum_{\tau\vdash m}\frac{\omega_\lambda(\tau)\omega_\mu(\tau)\chi_j(C_{2\tau})}{z^\tau},
\end{align*}
where $\omega_\lambda$ is the character of $S_m$ representation $L_{\lambda}$ and $C_{2\tau}$ an element of cycle type $2\tau =(2\tau_1, ..., 2\tau_l)$.
If $j$ is odd, then $\chi_j(C_{2\tau})=0$. When $j=2r$, we have $\chi_j(C_{2\tau})=(-1)^rE_t(\tau)$, and $E_t$ is the character of the wedge module $\wedge^t\mbb C^m$.
Then
$$\Delta_\lambda \otimes \chi_j=\sum_{\tau \vdash r, \nu \vdash m-r}(-1)^rc_{\nu, \tau}^\lambda c_{\nu, \tau^t}^\mu \Delta_{\mu}.$$
The multiplicity of $L_\mu$ as an irreducible constituent in $L_\lambda\otimes \wedge^t \mbb C^m$ is $$\sum_{\tau \vdash r, \nu \vdash m-r}c_{\nu, \tau}^\lambda c_{\nu, \tau^t}^\mu.$$

On the other hand,  we have
 $\chi_{\lambda^\pm}=\tfrac{1}{2}({\chi_\lambda \pm \Delta_\lambda})$. Combining this with the results of Theorem \ref{theorem 1.1.2} and Theorem \ref{bipartition multiplication of B} and , we obtain the desired formula.
\end{proof}

\section{The exceptional cases}
We now consider the types $E_6$, $E_7$ and $E_8$, and adopt the same notation as in \cite[\S 4]{BM99}. All results are obtained through computations in \emph{SageMath} \cite{S}. The simple roots are chosen as follows (for $E_6$ we take the first $6$, for $E_7$ the first $7$, and for $E_8$ all $8$):
\begin{align*}
	&\alpha_1=\tfrac{1}{2}(\epsilon_1+\epsilon_8)-\tfrac{1}{2}(\epsilon_2+\epsilon_3+\epsilon_4+\epsilon_5+\epsilon_6+\epsilon_7), \\
	&\alpha_2=\epsilon_1+\epsilon_2, \ \ \alpha_3=\epsilon_2-\epsilon_1, \ \ \alpha_4=\epsilon_3-\epsilon_2, \\
	&\alpha_5=\epsilon_4-\epsilon_3, \ \ \alpha_6=\epsilon_5-\epsilon_4, \ \ \alpha_7=\epsilon_6-\epsilon_5, \ \ \alpha_8=\epsilon_7-\epsilon_6.
\end{align*}
Then $s_i(v)=v-2\frac{(v,\alpha_i)}{(\alpha_i,\alpha_i)}\alpha_i$ defines the simple reflections of the Weyl group. Define
$$ w_1=s_1, \ \ w_2=s_1s_2,  \ \ w_3 =s_1s_3, \ \ w_4=s_1s_2s_3s_1s_4s_2s_3s_4s_5s_4, \ \ w_6=s_2s_3s_4s_5.$$

 Throughout this section, we record the character values of the restriction to $W$ of each one $\rm K$-type representation,

\textbf{Type $E_6$.} There are seven one $\rm K$-type representations for type $E_6$.
The character values on the elements $1, w_1, w_3, w_2, w_4, w_6$
are as follows.
\begin{align*}
	\begin{array}{cccccc}
		1 &  w_1 & w_3 & w_2 & w_4 & w_6  \\
		1^{\pm} & \pm1 & 1 & 1 & 1 & 1 \\
		10  & 0 & -2 & 2 & 2 & 0 \\
		15^{\pm} & \pm5 & 0 & 3 & -1 & -2 \\
		24^{\pm} & \pm 4 & 0 & 0 & 0 & 2
	\end{array}
\end{align*}
A sign $\pm$ attached to a dimension indicates that the corresponding restriction to $W$ is not isomorphic to its twist by the sign character.
The characters of the one $\rm K$-type representations that occur in
$\chi\otimes \wedge^iV$ are as follows.
\begin{align*}
	\begin{array}{c|c|c|c}
		\hline
		& \chi \otimes \wedge^1V & \chi \otimes \wedge^2V & \chi \otimes \wedge^3V  \\
		\hline
		\chi_{1^+} & \emptyset & \emptyset & \emptyset  \\
		\hline
		\chi_{10} & \emptyset & \emptyset & \chi_{10}  \\
		\hline
		\chi_{15^+} & \emptyset & \emptyset & \chi_{24^+}+\chi_{24^-} \\
		\hline
		\chi_{24^+} & \emptyset & \chi_{24^-} & \chi_{15^+}+\chi_{15^-}  \\
		\hline
	\end{array}
\end{align*}
We record only the constituents that are restrictions of one $K$-type representations.

\textbf{Type $E_7$.} There are $18$ one $K$-type representations of type $E_7$. The character values on the elements $1, w_1, w_3, w_2, w_4, w_6$ are as follows.
\begin{align*}
	\begin{array}{cccccc}
		1 & w_1 & w_3 & w_2 & w_4 & w_6 \\
		1^{\pm} & \pm1 & 1 & 1 & 1 & 1 \\
		15^{\pm} & \pm 5 & 0 & 3 & -1 & -2  \\
		21^{\pm} & \pm11 &6 & 5 & 1 & 2  \\
		35^{\pm} & \pm15 & 5 & 7 & -1 & -1 \\
		70^{\pm} & \pm10 &-5 & 6 & 2 & -1  \\
		84^{\pm} & \pm 4 & -6 & 4 & 4 & 2 \\
		105^{\pm} & \pm 5 & 0 & -3 & -3 & 2  \\
		105_\ast^{\pm} & \pm25 & 0 & 9 & -3 & -4 \\
		210^{\pm} & \pm 10 &-15 & 10 & 6 & 1
	\end{array}
\end{align*}
There are two non-isomorphic one $\rm K$-type representations of dimension $105$, and twisting by the sign character does not interchange them.
The one $\rm K$-type constituents of $\chi\otimes\wedge^i V$ are as follows:
\begin{align*}
	\begin{array}{c|c|c|c}
		\hline
		& \chi \otimes \wedge^1V & \chi \otimes \wedge^2V & \chi \otimes \wedge^3V  \\
		\hline
		\chi_{1^{+}} & \emptyset & \emptyset & \emptyset  \\
		\hline
		\chi_{15^{+}} & \chi_{105_\ast^+} & \emptyset & \chi_{105^+}  \\
		\hline
		\chi_{21^{+}} & \emptyset & \emptyset & \emptyset \\
		\hline
		\chi_{35^{+}} & \emptyset & \emptyset & \emptyset  \\
		\hline
		\chi_{70^{+}} & \chi_{210^+} & \chi_{84^-} & \chi_{210^+}  \\
		\hline
		\chi_{84^{+}} & \chi_{210^-} & \chi_{70^-} & \chi_{84^-}+\chi_{210^-}  \\
		\hline
		\chi_{105^{+}} & \emptyset & \chi_{105_\ast^+} & \chi_{15^+}+\chi_{105^-}+\chi_{105_\ast^-}  \\
		\hline
		\chi_{105_\ast^{+}} & \chi_{15^+} & \chi_{105^+}+\chi_{105_\ast^+} & \chi_{105^-}  \\
		\hline
		\chi_{210^{+}} & \chi_{70^+}+\chi_{84^-}+\chi_{210^-} & \chi_{84^+}+2\chi_{210^+} & \chi_{70^+}+\chi_{84^-}+2\chi_{210^-}  \\
		\hline
	\end{array}
\end{align*}

\textbf{Type $E_8$.} There are $18$ one $\rm K$-type representations for type $E_8$. Sixteen of them form eight sign-twist pairs of dimensions $1, 50, 84, 175, 448, 525, 700, 972$, and the remaining two, of dimensions $168$ and $420$, are self-twisted, that is, fixed by the sign twist.
\begin{align*}
	\begin{array}{cccccc}
		1 & w_1 & w_3 & w_2 & w_4 & w_6 \\
		1^{\pm} & \pm1 & 1 & 1 & 1 & 1 \\
		50^{\pm} & \pm 20 & 5 & 10 & -2 & -3  \\
		84^{\pm} & \pm 42 & 21 & 20 & 4 & 5  \\
		175^{\pm} & \pm 35 & -5 & 15 & -1 & -5 \\
		448^{\pm} & \pm 112 & 16 & 32 & 0 & 0  \\
		525^{\pm} & \pm 105 & 30 & 5 & -7 & 6 \\
		700^{\pm} & \pm 70 & -20 & 20 & 0 & -4  \\
		972^{\pm} & \pm 162 & 0 & 36 & 0 & 0 \\
		168 & 0 & -12 & 8 & 8 & 4 \\
		420 & 0 & -30 & 20 & 12 & 2
	\end{array}
\end{align*}
The characters of the one $\rm K$-type constituents of $\chi \otimes \wedge^i V$ are as follows.
\begin{align*}
	\begin{array}{c|c|c|c|c}
		\hline
		& \chi \otimes \wedge^1V & \chi \otimes \wedge^2V & \chi \otimes \wedge^3V  & \chi \otimes \wedge^4V \\
		\hline
		\chi_{1^{+}} & \emptyset & \emptyset & \emptyset & \emptyset \\
		\hline
		\chi_{50^{+}} & \emptyset & \emptyset & \emptyset & \emptyset  \\
		\hline
		\chi_{84^{+}} & \emptyset & \emptyset & \emptyset & \emptyset \\
		\hline
		\chi_{175^{+}} & \emptyset & \chi_{700^+} & \emptyset & \emptyset  \\
		\hline
		\chi_{448^{+}} & \emptyset & \emptyset & \emptyset & \emptyset  \\
		\hline
		\chi_{525^{+}} & \emptyset & \emptyset & \emptyset & \chi_{525^+}+\chi_{525^-}  \\
		\hline
		\chi_{700^{+}} & \emptyset & \chi_{175^+} & \emptyset & \chi_{175^+}+\chi_{175^-}  \\
		\hline
		\chi_{972^{+}} & \emptyset & \emptyset & \emptyset & \chi_{972^+}+\chi_{972_\ast^+}   \\
		\hline
		\chi_{168} & \emptyset & \chi_{420} & \emptyset &  \chi_{168}+\chi_{420}  \\
		\hline
		\chi_{420} & \emptyset & \chi_{168}+\chi_{420} & \emptyset &  \chi_{168}+\chi_{420}  \\
		\hline
	\end{array}
\end{align*}

\textbf{Type $F_4$.} For the root system of type $F_4$ with simple roots
\begin{align*}
	\alpha_1=\epsilon_2-\epsilon_3,\quad \alpha_2=\epsilon_3-\epsilon_4,\quad \alpha_3=\epsilon_4,\quad \alpha_4=\tfrac{1}{2}(\epsilon_1-\epsilon_2-\epsilon_3-\epsilon_4),
\end{align*}
we denote by $s_i$ the simple reflection associated with $\alpha_i$. Take $h=k_{\alpha_1}$ and $k_{\alpha_3}=1$, the parameters being constant on each of the two root length orbits. Set
\begin{align*}
	w_3 &= s_1s_2,  \quad w_{3'} = s_3s_4, \\
	w_{12} &= s_1s_2s_3s_4, \\
	w_6 &= w_{12}^2, \quad w_{6'} = w_6s_1s_2s_4s_3,  \quad w_{6''} = w_{12}^3 s_1s_2, \\
	w_2 &= s_1s_3,  \quad w_{2'} = s_1s_2s_4s_3s_2s_1s_3s_4, \\
	w_4 &= s_2s_3, \quad w_{4'} = w_{12}^3.
\end{align*}
For all values of $h$, there are nine one $\rm K$-type representations. The character values on the elements are as follows.
\begin{align*}
	\begin{array}{ccccccccccccc}
		1 & s_1 & s_3 & w_2 & w_{2'} & w_3 & w_{3'} & w_4 & w_{4'} & w_6 & w_{6'}  &  w_{6''} &  w_{12} \\
		1^{\pm} & \pm 1 & \pm1 & 1 & 1 & 1 & 1 & 1 & 1  & 1 & 1 & 1 & 1\\
		1_a^{\pm} & \mp 1 & \pm1 & -1 & 1 & 1 & 1 & -1  & 1 & 1 & 1 & 1 & 1\\
		2^{\pm} & 0 & \mp 2 & 0 & 2 & 2 & -1 & 0  & 2 & -1 & -1 & 2 & -1 \\
		2_a^{\pm} & \mp 2 & 0 & 0 & 2 & -1 & 2 & 0 & 2  & -1 & 2 & -1 & -1\\
		4 & 0 & 0 & 0 & 4 & -2 & -2 & 0 & 4  & 1 & -2 & -2 & 1 \\
	\end{array}
\end{align*}
When $h=1$, there are five additional one $\rm K$-type representations, of dimensions $4, 4, 6, 8, 8$. Their character values are as follows.
\begin{align*}
	\begin{array}{ccccccccccccc}
		1 & s_1 & s_{3} & w_2 & w_{2'} & w_3 & w_{3'} & w_4 & w_{4'} & w_6 & w_{6'}  &  w_{6''} &  w_{12} \\
		4_a^{\pm} & \pm2 & \mp 2 & 0 & 0 & 1 & 1 & -2 & 0  & 2 & -1 & -1 & 0 \\
		6 & 0 & 0 &  2 & -2 & 0 & 0 & -2 & 2  & 3 & 0 & 0 & -1  \\
		8^{\pm} & \pm 4 & 0 & 0 & 0 & -1 & 2 & 0 & 0  & -2 & -2 & 1 & 0 \\
	\end{array}
\end{align*}
When $h=2$, the representations $4_a^{\pm}$ and $6$ are again one $\rm K$-type representations, while the two $8$-dimensional ones are not. Instead, two other \(8\)-dimensional representations, denoted \(8_a^{\pm}\), again differing from each other by the sign character, are one \(K\)-type representations. Their character values on the same elements are as follows.
\begin{align*}
	\begin{array}{ccccccccccccc}
		1 & s_1 & s_3 & w_2 & w_{2'} & w_3 & w_{3'} & w_4 & w_{4'} & w_6 & w_{6'}  &  w_{6''} &  w_{12} \\
		8_a^{\pm} & 0 & \pm4 & 0 & 0 & 2 & -1 & 0 & 0  & -2 & 1 & -2 & 0 \\
	\end{array}
\end{align*}
The one $K$-type constituents of $\chi\otimes\wedge^i V$ are as follows:
\begin{align*}
	\begin{array}{c|c|c}
		\hline
		& \chi \otimes \wedge^1V & \chi \otimes \wedge^2V   \\
		\hline
		\chi_{1^{+}} & \emptyset & \emptyset   \\
		\hline
		\chi_{1_a^{+}} & \chi_{4_a^-} & \chi_{6}   \\
		\hline
		\chi_{2^{+}} & \chi_{8_a^-} & \emptyset  \\
		\hline
		\chi_{2_a^{+}} & \chi_{8^-} & \emptyset   \\
		\hline
		\chi_{4} & \emptyset & \chi_{6}   \\
		\hline
		\chi_{4_a^{+}} & \chi_{1_a^-}+\chi_6 & \chi_{4_a^+}+\chi_{4_a^-}   \\
		\hline
		\chi_{6} & \chi_{4_a^+}+\chi_{4_a^-}  & \chi_{1_a^+}+\chi_{1_a^-}+\chi_{4}+\chi_{6}   \\
		\hline
		\chi_{8^+} & \chi_{2_a^-} & \chi_{8^+}+\chi_{8^-}+\chi_{8_a^+}+\chi_{8_a^-}   \\
		\hline
		\chi_{8_a^{+}} & \chi_{2_a^-} & \chi_{8^+}+\chi_{8^-}+\chi_{8_a^+}+\chi_{8_a^-}  \\
		\hline
	\end{array}
\end{align*}

\textbf{Type $G_2$.} Consider the root system of type $G_2$ realized in the hyperplane
\begin{align*}
	H=\{x \in \mbb R^3: x_1+x_2+x_3=0 \},
\end{align*}
with simple roots
\begin{align*}
	\alpha_1=\varepsilon_1-\varepsilon_2 \ ({\rm short}),  \ \ \ \alpha_2=-\varepsilon_1+2\varepsilon_2-\varepsilon_3 \ ({\rm long}).
\end{align*}
Let $s_i$ be the simple reflection associated with $\alpha_i$. We set $k_{\alpha_1}=1$ for the short simple root and write $h=k_{\alpha_2}$ for the long simple root.
 Define $c=s_1s_2$.
For all values of $h$, there are four one $\rm K$-type representations are one-dimensional, and when $h=1$ or $h=3$ there is an additional one $K$-type representation of dimension two which is not the natural representation $V$.
 The character values of these are listed as follows.
\begin{align*}
	\begin{array}{cccccc}
		1 & s_1 & s_2 & c & c^2 & c^3  \\
		1^{\pm} & \pm 1 & \pm1 & 1 & 1 & 1 \\
		1_\varepsilon^{\pm} & \pm 1 & \mp 1 & -1 & 1 & -1 \\
		2 & 0 & 0 & -1 & -1 & 2
	\end{array}
\end{align*}
The characters of the one $\rm K$-type constituents of $\chi\otimes V$ are as follows.
\begin{align*}
	\begin{array}{c|c}
		\hline
		& \chi \otimes V   \\
		\hline
		\chi_{1^{+}} & \emptyset   \\
		\hline
		\chi_{1_\varepsilon^{+}} & \chi_{2}    \\
		\hline
		\chi_{2} & \chi_{1_\varepsilon^+}+\chi_{1_\varepsilon^-}  \\
		\hline
	\end{array}
\end{align*}


\begin{thebibliography}{99999}\frenchspacing
\bibitem[Bo76]{Bo76}
A. Borel, \emph{Admissible representations of a semi-simple group over a local field with vectors fixed under an Iwahori subgroup}, Invent. Math. 35 (1976), 233-259. \url{https://doi.org/10.1007/BF01390139}

\bibitem[BM89]{BM89}
D. Barbasch, A. Moy, \emph{A unitarity criterion for $p$-adic groups}, Invent. Math. 98 (1989),  19-3. \url{https://doi.org/10.1007/BF01388842}

\bibitem[BM93]{BM93}
D. Barbasch, A. Moy, \emph{Reduction to real infinitesimal character in affine Hecke algebras}, J. Amer. Math. Soc. 6 (1993), 611-635. \url{https://doi.org/10.1090/S0894-0347-1993-1186959-0}

\bibitem[BM99]{BM99}	
D, Barbasch, A. Moy, \emph{Classification of one $\rm K$-type representations}, Trans. Amer. Math. Soc. 351 (1999), 4245-4261. \url{https://doi.org/10.1090/S0002-9947-99-02171-6}

\bibitem[BW00]{BW00}
A. Borel, N. Wallach, \emph{Continuous cohomology, discrete subgroups, and representations of reductive groups}, Second edition, Math. Surveys Monogr. 67, Amer. Math. Soc., Providence, RI, 2000.
\url{https://doi.org/10.1090/surv/067}

\bibitem[Ca74]{Ca74}
W. Casselman, \emph{On a $p$-adic vanishing theorem of Garland}, Bull. Amer. Math. Soc. 80 (1974), 1001-1004. \url{https://doi.org/10.1090/S0002-9904-1974-13611-6}

\bibitem[Ch16]{Ch16}
K.Y. Chan, \emph{Duality for Ext-groups and extensions of discrete series for graded Hecke algebras}, Adv. Math, 294 (2016), 410-453. \url{http://dx.doi.org/10.1016/j.aim.2016.03.002}

\bibitem[Fu97]{Fu97}
W. Fulton, \emph{Young Tableaux, with Applications to Representation Theory and Geometry}, London Mathematical Society Student Texts, vol 35, Cambridge University Press, Cambridge, 1997.
\url{https://doi.org/10.1017/CBO9780511626241}

\bibitem[Ha15]{Ha15}
   T. Hayashi, \emph{A decomposition rule for certain tensor product
representations of the symmetric groups}, J. Algebra 434 (2015),  46每64. \url{http://dx.doi.org/10.1016/j.jalgebra.2015.03.016}

\bibitem[HO97]{HO97}
G. J. Heckman, E. M. Opdam, \emph{Yang's system of particles and Hecke algebras}, Ann. of Math. (2) 145 (1997), 139-173. \url{https://doi.org/10.2307/2951825}

\bibitem[KL87]{KL87}
D. Kazhdan, G. Lusztig, \emph{Proof of the Deligne-Langlands conjecture for Hecke algebras}, Invent. Math. 87 (1987), no. 1, 153-215. \url{https://doi.org/10.1007/BF01389157}

\bibitem[Lu79]{Lu79}
G. Lusztig, \emph{A class of irreducible representations of the Weyl group}, Nederl. Akad. Wetensch. Indag. Math. 82 (1979), 323-335. \url{https://doi.org/10.1016/1385-7258(79)90036-2}

\bibitem[Lu88]{Lu88}
G. Lusztig, \emph{Cuspidal local systems and graded Hecke algebras, I}, Publ. Math. IHES 67 (1988), 145-202. \url{https://doi.org/10.1007/BF02699129}

\bibitem[Lu89]{Lu89}
G. Lusztig, \emph{Affine Hecke algebras and their graded versions}, J. Amer. Math. Soc. 2 (1989), 599每635. \url{https://doi.org/10.2307/1990945}

\bibitem[Lu95]{Lu95}
G. Lusztig, \emph{Cuspidal local systems and graded Hecke algerbas. II} in \emph{Representations of groups}, (B. Allison and G. Cliff, eds.), CMS Conf. Proc., 16, Amer. Math. Soc., Providence, RI, 1995, pp217-275.

\bibitem[Ok98]{Ok98}
S. Okada, \emph{ Applications of Minor Summation Formulas to Rectangular-Shaped Representations of Classical Groups}, J. Algebra 205 (1998), 337-367. \url{https://doi.org/10.1006/jabr.1997.7408}

\bibitem[Op95]{Op95}
E. M. Opdam, \emph{Harmonic analysis for certain representations of graded Hecke algebras}, Acta Math. 175 (1995), 75-121. \url{https://doi.org/10.1007/BF02392487}
	
\bibitem[Or05]{Or05}
S. Orlik, \emph{On extensions of generalized Steinberg representations}, J. Algebra 293 (2005), 611-630. \url{https://doi.org/10.1016/j.jalgebra.2005.03.028}

\bibitem[OS10]{OS10}
E. M. Opdam, M. Solleveld, \emph{Discrete series characters for affine Hecke algebras and their formal degrees}, Acta Math. 205 (2010), 105-187. \url{https://doi.org/10.1007/s11511-010-0052-9}

\bibitem[S]{S}
The Sage Developers, \emph{SageMath, the Sage Mathematics Software System}, \url{https://www.sagemath.org}.

\bibitem[St12]{St12}
R. P. Stanley, \emph{Enumerative Combinatorics, Vol. 1}, 2nd ed.,
Cambridge University Press, Cambridge, 2012. \url{https://doi.org/10.1017/CBO9781139058520}
	
\bibitem[Ta15]{Ta15}
J. Taylor, \emph{Induced characters of type $D$ Weyl groups and the Littlewood每Richardson rule}, J. Pure Appl. Algebra 219 (2015), 3445-3452. \url{https://doi.org/10.1016/j.jpaa.2014.12.004}

\bibitem[VZ84]{VZ84}
D. A. Vogan, Jr., G. J. Zuckerman, \emph{Unitary representations with non-zero cohomology}, Compositio Math. 53 (1984), 51-90. \url{http://archive.numdam.org/item/CM_1984__53_1_51_0/}
	\end{thebibliography}
\end{document}